\documentclass[USenglish]{article}	
\usepackage[utf8]{inputenc}				%(only for the pdftex engine)
\usepackage[big,online]{dgruyter}	%values: small,big | online,print,work
\usepackage{lmodern} 
\usepackage{microtype}
\usepackage[numbers,square,sort&compress]{natbib}
\usepackage[dvipsnames]{xcolor}

\theoremstyle{dgthm}
\newtheorem{theorem}{Theorem}
\newtheorem{corollary}{Corollary}
\newtheorem{Pro}{Proposition}
\newtheorem{lemma}{Lemma}

\theoremstyle{dgdef}
\newtheorem{definition}{Definition}

\newtheorem{remark}{Remark}

\begin{document}

%%%--------------------------------------------%%%
	\articletype{Research Article}
	\received{Month	DD, YYYY}
	\revised{Month	DD, YYYY}
  \accepted{Month	DD, YYYY}
  \journalname{De~Gruyter~Journal}
  \journalyear{YYYY}
  \journalvolume{XX}
  \journalissue{X}
  \startpage{1}
  \aop
  \DOI{10.1515/sample-YYYY-XXXX}
%%%--------------------------------------------%%%

\title{Inverse Problems of Identifying the Unknown Transverse Shear Force in the Euler-Bernoulli Beam with Kelvin-Voigt Damping}
\runningtitle{Inverse Problems of Identifying the Unknown Transverse Shear Force }
%\subtitle{Insert subtitle if needed}

\author*[1]{K. Sakthivel}
%\ use * to mark the author as the corresponding author
\author[2]{A. Hasanov}
\author[3]{D. Anjuna} 
\runningauthor{ K. Sakthivel et al.}
\affil[1]{\protect\raggedright 
Indian Institute of Space Science and Technology, Department of Mathematics, Trivandrum- 695 547, India, e-mail:sakthivel@iist.ac.in}
\affil[2]{\protect\raggedright 
 Department of Mathematics, Kocaeli University, Turkey\\
 $~$ \c{S}ehit Ekrem Mah., Altin\c{s}ehir Sk., Ayazma Villalari, 22,
 Bah\v{c}ecik, Kocaeli - 41030, Turkey, 
 e-mail:  alemdar.hasanoglu@gmail.com }
	\affil[3]{\protect\raggedright
		 Indian Institute of Space Science and Technology, Department of Mathematics, Trivandrum- 695 547, India,  e-mail:  dileepanjuna10@gmail.com}
%\communicated{...}
%\dedication{...}
	
\abstract{In this paper, we study the inverse problems of determining the unknown transverse shear force $g(t)$ in a system governed by the damped Euler-Bernoulli equation $\rho(x)u_{tt}+\mu(x)u_t+ (r(x)u_{xx})_{xx}+ (\kappa(x)u_{xxt})_{xx}=0, ~(x,t)\in (0,\ell)\times(0,T],$ subject to the boundary conditions $u(0,t) =0$, $u_{x}(0,t)=0$, $\left[r(x)u_{xx}+\kappa(x)u_{xxt}\right]_{x=\ell}  =0$, $-\left[\big(r(x)u_{xx}+\kappa(x)u_{xxt}\big)_{x}\right]_{x=\ell}=g(t)$, $t\in [0,T]$, from the measured deflection $\nu(t):=u(\ell,t)$, $t \in [0,T]$, and from the bending moment $\omega(t):=-\left( r(0)u_{xx}(0,t)+\kappa(0)u_{xxt}(0,t) \right)$, $t \in [0,T]$, where the terms $(\kappa(x)u_{xxt})_{xx}$ and $\mu(x)u_t$ account for the Kelvin-Voigt damping and external damping, respectively.
	The main purpose of this study is to analyze the Kelvin-Voigt damping effect on determining the unknown transverse shear force (boundary input) through the given boundary measurements. The inverse problems are transformed into minimization problems for Tikhonov functionals, and it is shown that the regularized functionals admit unique solutions for the inverse problems. By suitable regularity on the admissible class of shear force $g(t),$ we prove that these functionals are Fr\'echet differentiable, and the derivatives are expressed through the solutions of corresponding adjoint problems posed with measured data as boundary data associated with the direct problem. The solvability of these adjoint problems is obtained under the minimal regularity of the boundary data $g(t)$, which turns out to be the regularizing effect of the Kelvin-Voigt damping in the direct problem. Furthermore, using the Fr\'echet derivative of the more regularized Tikhonov functionals,   we obtain remarkable Lipschitz stability estimates for the transverse shear force in terms of the given measurement by a feasible condition only on the Kelvin-Voigt damping coefficient.}


\maketitle

\section{Introduction} 
\label{intro}
A beam is a structural element made to resist forces acting laterally to its axis. Daniel Bernoulli derived the equation of motion for transverse vibration of thin beams in 1735, while Euler provided the first solutions for various support conditions in 1744. The Euler-Bernoulli beam theory is a simplified version of linear elasticity theory that describes the relationship between deflection and applied load. Many mechanical systems from industry and engineering use the Euler–Bernoulli beam equation to represent bending vibration \cite{LMeirovitch}. Analysis and simulation of such systems have become key research areas because of the necessity to manage the dynamics of these systems. The static equation $ (EI u^{\prime \prime}(x))^{\prime \prime}=q(x),$ where $ u$ is the beam deflection, $EI$ is the flexural rigidity and $ q $ is the distributed load (see, \cite{timoshenkostability},\cite{rao2007}), is a basic model for the Euler-Bernoulli beam equation subject to external load.

In the past two decades, several studies have been conducted on inverse problems corresponding to this classical model with appropriate boundary conditions  arising from a wide variety of basic applications in civil, mechanical, and aeronautical engineering (see,  \cite{cheng2006innermass}, \cite{fryba1996railway}, \cite{fryba1999vibration}, \cite{vinod2007freevibration}). The classical Euler-Bernoulli beam model with appropriate damping mechanisms such as viscous (air) damping, strain rate damping, spatial hysteresis, and time hysteresis play a significant role in applications (see,  \cite{banks1991}). A general  model of the damped beam equation with generic damping is given by $ u_{tt}(x,t)+L_1 u_{t}(x,t)+L_2 u(x,t)+\left((EI(x)/\rho) u_{xx}\right)_{xx}=q(x,t),$ where the term
$L_1u_t(x,t) + L_2u(x,t)$ accounts for the damping mechanisms of this model. External damping
mechanisms usually determine the nature of coefficient $L_1 $ while the internal damping mechanisms often determine the coefficient $L_2. $ While the inverse problems related to  the
classical (undamped) beam model with appropriate boundary conditions have been well studied over the past two decades, a very limited study has only been done in the damped beam equation. It should be noted here that the unique determination of transverse load or spatial load from the final time displacement for undamped wave or beam model was not feasible (see, \cite{hasanouglu2017introduction}). Recently we proved that (see, \cite{hasanov2021},\cite{anjuna:2021}) in the presence of the external damping, though we didn't consider the internal damping in these papers, we can uniquely recover source functions from the final time measured data under suitable conditions of the temporal load. It is understood that the nature of the damping term drastically changes the nature of the solutions to the direct beam model, which in turn helps to obtain a feasible solution to the inverse problems as well. By taking this motivating factor into account, in this paper, we analyze the inverse problem of determining the unknown transverse shear force in the presence of both internal and external damping factors. It is worth noting that the wave equation with Kelvin-Voigt damping in a bounded domain was studied in \cite{ammari2020}. 

Consider the Euler-Bernoulli beam with more general physical coefficients and damping effects as follows  (see, \cite{Clough:Penzien},  Chapter 17,  Section 4 and also refer to \cite{banks1991}):
\begin{eqnarray}\label{1-1}
	\left\{ \begin{array}{lc}
		\rho(x)u_{tt}+\mu(x)u_t+ (r(x)u_{xx})_{xx}+ (\kappa(x)u_{xxt})_{xx}=0, ~(x,t)\in \Omega_T, \\ [2pt]
		u(x,0)=0,~u_t(x,0)=0, ~  \ \  x\in (0,\ell),\\
		u(0,t) \ =0,~  u_{x}(0,t)=0, ~ \ \  t\in [0,T],\\ [2pt]
		\left[r(x)u_{xx}+\kappa(x)u_{xxt}\right]_{x=\ell} \; =0, \\ 
		\quad\quad		 -\left[\big(r(x)u_{xx}+\kappa(x)u_{xxt}\big)_{x}\right]_{x=\ell}=g(t),~  t\in [0,T],
	\end{array}\right.
\end{eqnarray}
where $\Omega_T:=(0,\ell)\times(0,T]$, $r(x):=E(x)I(x)>0$ is the flexural rigidity (or bending stiffness) of a non homogeneous beam while $E(x)>0$ is the elasticity modulus and $I(x)>0$ is the moment of inertia. The coefficient $\kappa(x):=c_dI(x)$ represents energy dissipated by friction internal to the beam,  where $c_d$  is the strain-rate damping coefficient.  The nature of the terms $\mu(x)u_t$ and $(\kappa(x)u_{xxt})_{xx}$ are determined by external and internal damping mechanisms, respectively. The non-negative coefficient $\mu(x)$ and the positive coefficient $\kappa(x)$ are called the viscous external damping and the strain-rate or Kelvin-Voigt damping coefficients, respectively.  For this considered model, there can be only two types of force effects acting on the right end $x=\ell,$ one is bending moment, and the other is transverse shear force. In this study, we look at the case in which the vibration is caused by an unknown transverse shear force $g(t):= -\left[\left(r(x)u_{xx}+\kappa(x)u_{xxt}\right)_x\right]$ at $x=\ell,$ which needs to be identified from either given measured deflection at  $x = \ell$ or bending moment at $ x=0.$
\par The following two inverse boundary-value problems (IBVPs) are formulated to the model (\ref{1-1}).

\noindent \textbf{IBVP-1.} \emph{Find the unknown transverse shear force $g(t)$ from measured deflection $\nu(t)$  given at the right end of the beam $x=\ell$:}
\begin{eqnarray}\label{1-2}
	\nu(t):=u(\ell,t), ~ t \in [0,T].
\end{eqnarray}
For a given $ g(t)$ from the set of admissible transverse shear forces, the problem (\ref{1-1}) with the solution $ u(x,t):=u(x,t;g)$ is  known as the direct problem. Here, the following tricky question arises. Can we swap the boundary conditions  $-\left[\big(\kappa(x)u_{xxt}+r(x)u_{xx}\big)_{x}\right]_{x=\ell}=g(t),$ $\nu(t):=u(\ell,t)$, and solve the initial boundary value problem
\begin{eqnarray}\label{1-3}
	\left\{ \begin{array}{lc}
		\rho(x)u_{tt}+\mu(x)u_t+(r(x)u_{xx})_{xx}+ (\kappa(x)u_{xx})_{xxt}=0, ~(x,t)\in \Omega_T, \\ [2pt]
		u(x,0)=0,~u_t(x,0)=0, ~ ~x\in (0,\ell),\\
		u(0,t) \; =0,~   u_{x}(0,t)=0, ~~ \ t\in [0,T],\\
		u(\ell,t)\;=\nu(t),~\left[r(x)u_{xx}+\kappa(x)u_{xxt}\right]_{x=\ell} =0, ~ t\in [0,T]
	\end{array}\right.
\end{eqnarray}
and then find the unknown transverse shear force $g(t)$ using the formula
$g(t)=-\left[\big(r(x)u_{xx}+\kappa(x)u_{xxt}\big)_{x}\right]_{x=\ell}\,?$
Therefore, a situation arises as if there is no need for any inverse problem. Besides, swapping of the above mentioned boundary conditions is also a mathematically correct approach.  However,  in this paper,  we consider the inverse  problem \eqref{1-1}-\eqref{1-2} as   the physical justification to that  model  comes from the fact that problem (\ref{1-1}) is a developed mathematical model of vibration of the cantilever tip due to the shear force interaction in Transverse Dynamic Force Microscope (TDFM) (see, \cite{antognozzi:2000}, \cite{antognozzi:binder:2001}, \cite{nguyen2001} and references therein).  The problem of determining the shear force is of great importance when specimen images and mechanical properties need to be computed at some submolecular precision (see, \cite{nguyen2001}).

The problem IBVP-1 defined by (\ref{1-1}) and (\ref{1-2}) can be reformulated as the invertibility of the Neumann-to-Dirichlet operator
\begin{eqnarray}\label{1-4}
	\left. \begin{array}{lc}
		\Phi :\mathcal{G}_1 \subset H^1(0,T)\mapsto L^2(0,T),  \ \ (\Phi g)(t):=u(\ell,t;g), ~ t \in [0,T], \\ [2pt]
		\mathcal{G}_1=\{ g\in H^1(0,T): g(0)=0,\ \Vert g \Vert_{H^1(0,T)}\leq C_{g},\, C_g>0\},
	\end{array}\right.
\end{eqnarray}
where $ \mathcal{G}_1$ is called the set of admissible inputs (shear forces). With the help of noise free measured output $ \nu(t),$ we can reformulate IBVP-1 in terms of functional equation as \begin{eqnarray}\label{2.0}
	\Phi g(t)=\nu(t), \ \nu \in L^2(0,T).
\end{eqnarray}
\par We note that the exact equality in (\ref{2.0}) can hold only in the case of noiseless measured output $\nu(t)$. However, it should be emphasized that in practice the measured output $\nu(t)$ always contains measurement errors, and hence exact equality in the functional equation (\ref{2.0}) is not possible.  Therefore, we introduce the Tikhonov functional to solve the minimization problem
\begin{eqnarray}\label{26}
	\min_{g\in \mathcal{G}_1} \mathcal{J}_1(g), \ \ \mathcal{J}_1(g):= \frac{1}{2}\,\Vert \Phi g-\nu \Vert^2_{L^2(0,T)},
\end{eqnarray}
whose solution,  according to \cite{Ivanov},  is defined as a quasi-solution of the inverse problem. Then we consider the same problem for the regularized Tikhonov functional
\begin{eqnarray}\label{54}
	\mathcal{J}_{1\alpha}(g):=\frac{1}{2}\,\Vert \Phi g -\nu\Vert^2_{L^2(0,T)}+\frac{\alpha}{2}\,\Vert g' \Vert^2 _{L^2(0,T)},
\end{eqnarray}
where $ \alpha >0$ is the parameter of regularization.

The second inverse problem, we study in this paper is formulated as follows. \\
\noindent \textbf{IBVP-2.} \emph{Find the unknown transverse shear force $g(t)$ from measured bending moment $ \omega(t)$  given at the beginning of the beam $x=0$:}
\begin{eqnarray}\label{1-5}
	\omega(t):=-\left( r(0)u_{xx}(0,t)+\kappa(0)u_{xxt}(0,t) \right), ~ t \in [0,T].
\end{eqnarray}
The problem IBVP-2 defined by (\ref{1-1}) and (\ref{1-5}) is related to the invertibility of
Neumann-to-Neumann operator
\begin{eqnarray}\label{1-6}
	\left.\begin{array}{lc}
		\Psi :\mathcal{G}_3 \subset H^3(0,T)\mapsto L^2(0,T), \\ [2pt]    (\Psi g)(t):=- \left(r(0)u_{xx}(0,t;g)+\kappa(0)u_{xxt}(0,t;g)\right),  t \in [0,T], \\ [2pt]
		\mathcal{G}_3=\{g\in H^3(0,T):g(0)=g^\prime(0)= g^{\prime \prime}(0)=0,  \Vert g \Vert_{H^3(0,T)}\leq \tilde C_{g},  \tilde C_{g}>0 \}.
	\end{array}\right.
\end{eqnarray}
In terms of the functional equation, we will again express the inverse problem (\ref{1-1}) and (\ref{1-5}) as follows
\begin{eqnarray}\label{3.0}
	\Psi g(t)=\omega(t), \ \ \omega\in L^2(0,T).
\end{eqnarray}
\par As in the case of IBVP-1, when the measured data $ \omega(t)$ contain random noise the exact equality in (\ref{3.0}) is not feasible. In this case we solve the  minimization problem for the Tikhonov functional
\begin{eqnarray}\label{5.9}
	\min_{g\in \mathcal{G}_3} \mathcal{J}_2(g), \ \ \ \mathcal{J}_2(g):= \frac{1}{2}\,\Vert \Psi g-\omega \Vert^2_{L^2(0,T)}
\end{eqnarray}
and the regularized Tikhonov functional is considered  as follows
\begin{eqnarray}\label{55}
	\mathcal{J}_{2\alpha}(g):=\frac{1}{2}\,\Vert \Psi g -\omega\Vert^2_{L^2(0,T)}+\frac{\alpha}{2}\,\Vert g''' \Vert^2 _{L^2(0,T)}.
\end{eqnarray}

It is worth noting that  IBVP-2 can be formulated using the admissible source $\mathcal{G}_1$ and the regularized  Tikhonov functionals $ \mathcal{J}_{1\alpha},$ $ \mathcal{J}_{2\alpha}$ can be defined  with usual  $L^2$ norm regularizer $\Vert g \Vert^2 _{L^2(0,T)}.$ Further,   the solvability of these inverse problems (see, Theorem \ref{45} or Remark \ref{Re1}) do not require the more regularized functionals as in (\ref{54}) and (\ref{55}), while these regularized functionals  are crucial to derive the stability estimates (see, Theorems \ref{mnk}, \ref{st2}). We study the inverse problems (\ref{1-1})-(\ref{1-2}) and (\ref{1-1}),(\ref{1-5}) as a minimization problems for the Tikhonov functionals $ \mathcal J_{1\alpha}(g)$ and $\mathcal  J_{2\alpha}(g)$ on the set $ \mathcal{G}_1,$ $ \mathcal{G}_3$ respectively. In the absence of the internal damping term $(\kappa(x)u_{xxt})_{xx}$ in (\ref{1-1}), the inverse problems of (\ref{1-1}) with measurements (\ref{1-2}) and (\ref{1-5}) were studied respectively  in \cite{AH:OB:CS-2019} and \cite{AH:OB:HI-2019}.

As we mentioned earlier, the inverse problem corresponding to the classical Euler- Bernoulli beam equation is a vast field of study. Early studies attempted to identify the unknown coefficients $A(x)$ and $I(x)$ in the simplest Euler-Bernoulli equation $\rho A(x)u_{tt} + (EI(x)u_{xx})_{xx} = 0$ from spectral data \cite{gldwell1986spectral}. Later works studied inverse problems more relevant to applications by using additional realistic data; specifically, boundary measured output data, since computing spectral data is difficult in practice. In \cite{chang2007identification}, the authors established the uniqueness of the identification of spatial density $\rho(x)$ and inertia $r(x)$ from the boundary data for the cantilever Euler-Bernoulli beam. The identification of an unknown spatial load using final time-measured deflection $ u(x,T)$ or final time measured velocity $ u_{t}(x,T)$ was studied in \cite{hasanov2009identification}. In contrast, the temporal load in the beam equation is identified with the help of measured boundary data $u_{x}(0,t),$ that is measured slope at the end $ x=0$ in \cite{hasanov2016identification}, and numerical algorithm for finding unknown sources has been developed. Two inverse source problems of detecting asynchronously distributed spatial load in $\rho(x)u_{tt}+\mu(x)u_{t}+(k(x)u_{xx})_{xx}-T_{r}u_{xx}=\sum_{m=1}^{M}h_{m}(t)f_{m}(x) $ with hinged-clamped endpoints are investigated in the work \cite{hasanovkawano2016identification}. For further results on the inverse source problems of the Euler-Bernoulli beam and plate equations,  one may refer to \cite{M. Grimmonprez}, \cite{kawano},  \cite{nicaise2004},  \cite{K. Van Bockstal}.   

Next, let us review some of the recent papers on the Euler-Bernoulli equation with external damping. In the paper \cite{AH:OB:CS-2019}, the authors determines the unknown transverse shear force by using measured boundary deflection $ u(\ell,t)$, and in the article \cite{AH:OB:HI-2019}, they consider the same inverse problem based on  measured bending moment$ -r(0)u_{xx}(0,t).$ When the temporal load $G(t) = 1$ and $G(t) = e^{-\eta t},\eta>0$ the effect of the damping parameter in the unique determination of spatial load in the Euler-Bernoulli beam equation from final time measured data has been explored in \cite{hasanov2021} by applying the  singular value decomposition. The numerical reconstruction of this problem with temporal load $G(t)=\cos(\omega t)$ was investigated in the recent paper \cite{HRB}.  In addition to these literature, there are some classic papers on the inverse problem of the Euler-Bernoulli beam equation with the Kelvin-Voigt damping or viscous damping. In \cite {graif} and \cite{banks1990}, the parameter identification of the Euler-Bernoulli beam equation with structural or viscous damping was investigated, and the numerical approximations of those quantities were studied. The paper \cite{itokekelvin} determined the stiffness $ EI(x)$, damping coefficient $DI(x),$ and initial data of Euler Bernoulli beam equation $u_{tt}+\left(EI(x)u_{xx}+DI(x)u_{xxt}\right)_{xx}=f(x,t),\ t>0$ using the spectral data of the model problem.

In this paper, we generalize the existing mathematical model for the Euler-Bernoulli equation by including all the possible physical coefficients together with both damping effects given by external damping as well as Kelvin-Voigt damping.  In contrast to external damping, the inclusion of the Kelvin-Voigt damping affects the free-end boundary conditions depending on the moment as it is strain-dependent. The presence of Kelvin-Voigt damping and the mixed boundary conditions, in turn, makes the problem more complicated, and the boundary data (shear force) determination under these conditions becomes difficult.  Indeed, for the existence and uniqueness of the solution to the direct problem (\ref{1-1}), we need to develop appropriate identities to handle these types of mixed boundary conditions for deriving priori estimates for the direct problem and proving the regularity of solutions. This makes the current paper different from  \cite{AH:OB:CS-2019}, \cite{AH:OB:HI-2019} where only external damping $\mu(x) u_t$ is considered, and so the balancing effect of the Kelvin-Voigt damping term on the  boundary conditions is also not needed on those papers.  On the other hand, Kelvin-Voigt damping also has some sort of regularizing effect in proving the solutions of the direct problem. For instance, in \cite{AH:OB:CS-2019}, \cite{AH:OB:HI-2019}, the authors require higher regularity like $g\in H^2(0,T)$ to prove $ u_{t}\in L^2(0,T;\mathcal{V}_1^2(0,\ell)), u_{tt}\in L^2(0,T;L^2(0,\ell)),$ whereas in this paper we prove those estimates with $ g\in H^1(0,T)$ by coupling with appropriate identities and Sobolev embedding theorems.  We also prove the existence of solutions to the inverse problems (\ref{1-1})-(\ref{1-2}) and (\ref{1-1}),(\ref{1-5}) when the transverse shear force $g(t)$ belongs to the admissible inputs $\mathcal G_1 \subset H^1(0,T).$

The main contribution of the paper is summarized as follows:
\begin{itemize}
	\item The existence and uniqueness of the weak and regular weak solutions to the direct problem    are proved. Furthermore, the necessary a priori estimates are derived.
	\item Solvability of the inverse problems (\ref{1-1})-(\ref{1-2}) and (\ref{1-1}),(\ref{1-5}), defined as IBVP-1 and IBVP-2, and governed by Neuman-to-Dirichlet and Neumann-to-Neumann operators, respectively, are studied in appropriate admissible set of transverse shear sources $\mathcal G_m \subset H^m(0,T)$, $m=1,3$. It is demonstrated that for compactness and Lipschitz continuity, the operator $ \Psi$ does not require $ \mathcal{G}_3$ regularity. We just need $ \mathcal{G}_2$ regularity for $\Psi$'s compactness, and $ \mathcal{G}_1$ regularity for Lipschitz continuity, while these results can be verified  for the operator $ \Phi$ on the admissible source $\mathcal{G}_1$ itself.
	\item The Tikhonov functionals $\mathcal J_{m\alpha}(g)$, $m=1,2$ are introduced, and the Fr\'echet derivatives of these functionals are derived in through the solutions of corresponding adjoint problems. It is shown that for IBVP-1, the admissible source $g(t)$ needs to be in ${\mathcal G_1}$, while for IBVP-2 the more regular set of admissible sources ${\mathcal G_3}\subset  H^3(0,T)$ is needed (see, Section \ref{s1}). For IBVP-2, the regularity of the admissible sources is also needed  for the solvability of the  adjoint problem.
	\item Other remarkable results are the Lipschitz type stability estimates for IBVP-1 and IBVP-2. We provide a local stability estimate for the unknown shear force $ g \in \mathcal{G}_{1}$ when  Kelvin-Voigt damping coefficient $ \kappa(x)>0$ satisfies a condition on its lower bound, while the external damping coefficient is nonnegative: $ \mu(x)\geq 0.$ It should be noted that this stability result is valid even when the external damping effect is not present, that is, when $\mu(x)=0.$  In the case of IBVP-2,  we establish a stability estimate for $ g\in \mathcal{G}_3$ under a feasible condition on the parameter of regularization $\alpha$. Both the stability results are obtained when a more smooth regularization terms are added to the Tikhonov functionals $\mathcal J_{m}(g),m=1,2.$  These results give a new perspective for the stability analysis of the shear force determination in the presence of both damping terms.
\end{itemize}

The paper is organized as follows. The existence and uniqueness, and regularity of solutions are given in Section \ref{s3}. The solvability of inverse problems are given in Sections \ref{s5}. The Section \ref{s1} is devoted to the Fr\'echet derivative of the Tikhonov functionals and Lipschitz continuity of the Fr\'echet derivatives. The monotonocity of the gradient based algorithm is analyzed in Section \ref{s55}.   The stability analysis of the inverse source problems are discussed in Section \ref{s2}.
\section{Existence and Uniqueness  of Weak Solutions to  Direct Problem}\label{s3}
In this section, we consider problem (\ref{1-1}) with nonhomogeneous initial conditions:
\begin{eqnarray}\label{1-0}
	\left\{ \begin{array}{lc}
		\rho(x)u_{tt}+\mu(x)u_t+ (r(x)u_{xx})_{xx}+ (\kappa(x)u_{xxt})_{xx}=0, ~(x,t)\in \Omega_T, \\ [2pt]
		u(x,0)=u_0,~u_t(x,0)=v_0, ~ \qquad  x\in (0,\ell),\\ 
		u(0,t) \ =0,~ u_{x}(0,t)=0, ~ \qquad \; \quad  t\in [0,T],\\ [2pt]
		\left[r(x)u_{xx}+\kappa(x)u_{xxt}\right]_{x=\ell} \; =0, \\ 
		\quad \quad	-\left[\big(r(x)u_{xx}+\kappa(x)u_{xxt}\big)_{x}\right]_{x=\ell}=g(t), \  t\in [0,T].
	\end{array}\right.
\end{eqnarray}
In the study of direct and inverse problems related to the model (\ref{1-1}), the following basic  assumptions are used:
\begin{eqnarray} \label{3-5}
	\left \{ \begin{array}{ll}
		\rho, \mu \in L^\infty(0,\ell), \  r,\kappa \in H^2(0,\ell) , \ g\in H^1(0,T),\,  g(0)=0,\\ [2.5pt]
		0<\rho_0\leq \rho(x)\leq \rho_1,  \  \  \ 0<r_0\leq r(x)\leq r_1,\\ [2.5pt]
		0\leq \mu_0\leq \mu(x)\leq \mu_1, \ \ \ 0< \kappa_0 \leq \kappa(x)\leq \kappa_1.	\end{array} \right.
\end{eqnarray}
\begin{definition}\label{d1}
	Let $ 0<T<+\infty,$ $u_0 \in \mathcal{V}_1^2(0,\ell)$ and  $ v_0\in \mathcal{V}_1^2(0,\ell)$ be given.   We say a function $ u \in L^2(0,T;\mathcal{V}_1^2(0,\ell))$ with $ u_t \in L^2(0,T;\mathcal{V}_1^2(0,\ell))$ and  $ u_{tt} \in L^2(0,T;L^2(0,
	\ell))$ is a weak solution of (\ref{1-0}) provided
	\begin{eqnarray*}
		i)&& \ \ (\rho \ u_{tt}(t),v)+(\mu \ u_{t}(t),v)+(ru_{xx}(t),v_{xx}) \\  &&+(\kappa u_{xxt}(t),v_{xx})  =g(t)v(\ell), \ \forall \ v\in \mathcal{V}_1^{2}(0,\ell), \  \  \mbox{a.e.} \ t\in[0,T]. \\
		ii)&& \ \  u(0)=u_{0}, \ u_{t}(0)=v_{0},
	\end{eqnarray*}	
	where $ \mathcal{V}_1^2(0,\ell)=\{ v\in H^2(0,\ell):v(0)=v_{x}(0)=0\}$ and it's dual space is $ \mathcal{V}_1^2(0,\ell)^{\prime}$. The space $ \mathcal{V}^2_1(0,\ell)$ is equipped with the standard Sobolev norm
	$ \Vert v \Vert_{\mathcal{V}_1^2(0,\ell)}:= \Big ( \int_0^\ell \left(v^2+v_x^2+v_{xx}^2\right)dx\Big)^{\frac{1}{2}}.$
\end{definition}
For any $v\in\mathcal{V}_1^2(0,\ell),$ we obtain the Poincar\'e inequalities, $ \Vert v\Vert_{L^2(0,\ell)}\leq 2 \ell \Vert v_x\Vert _{L^2(0,\ell)}$ $\leq 4 \ell^2 \Vert v_{xx}\Vert_{L^2(0,\ell)}.$ Thus, we have
\begin{eqnarray}\label{4.8}
	\Vert v\Vert_{\mathcal{V}_1^2(0,\ell)}\leq \sqrt{C^{\ast}}\Vert v_{xx}\Vert_{L^2(0,\ell)},~C^{\ast}=4\ell^2(1+4\ell^2)+1.
\end{eqnarray}
It is clear that $ \Vert v\Vert_{\mathcal{V}_1^2(0,\ell)}$ is equivalent to $ \Vert v_{xx}\Vert_{L^2(0,\ell)}.$

From Definition \ref{d1}, it is evident that $ u\in H^1(0,T;\mathcal{V}_1^2(0,\ell)),$ so that $ u\in C([0,T];\mathcal{V}_1^2(0,\ell))$ and  $ u_t \in C([0,T];L^2(0,\ell)).$ Consequently,  the equalities $ u(0)=u_0$ and $ u_t(0)=v_0$ can be justified.

We apply the Fa\'{e}do-Galerkin approximation method to illustrate  that there exists a unique weak solution to direct problem (\ref{1-0}).
First, we choose a sequence of smooth functions $ \{\xi_i\}_{i=1}^{n},$ which form an orthonormal and orthogonal basis for $ L^2(0,\ell)$ and $ \mathcal{V}_1^2(0,\ell)$ respectively. Then, we construct the $ n$ dimensional subspace $ W_n:= span\{\xi_{1},\xi_{2},...,\xi_{n}\} $ of $ \mathcal{V}_1^{2}(0,\ell)$ and seek the Fa\'edo-Galerkin approximation $u_{n}(t):=u_{n}(x,t)$ of the form $ u_n(t)=\sum_{i=1}^{n} d_{i,n}(t)\xi_{i},$ $u_{0,{n}}=\sum_{i=1}^{n}p_{{i,{n}}}\xi_{i},$ and $v_{0,n}=\sum_{i=1}^{n}q_{i,{n}}\xi_{i},$
where we hope to find the coefficients $ d_{i,n},\ p_{i,n}$ and $ q_{i,n}$ so that
\begin{eqnarray} \label{w2}
	\left \{ \begin{array}{lc}
		(\rho \ u_{n}^{\prime \prime}(t),v)+(\mu \ u_{n}^{\prime}(t),v)+(r \ u_{n,{xx}}(t),v_{xx}) \\ +(\kappa \ u_{n,{xx}}^\prime(t),v_{xx})=g(t)v(\ell),\ \ \forall \ v\in W_{n}, \  \  t\in[0,T],  \\ [2pt]
		u_{n}(0)=u_{0,n}, \ \ u_{n}^{\prime}(0)=v_{0,n}.
	\end{array} \right.
\end{eqnarray}
By inserting $ v=\xi_{j}, \ j=1,2,3,...n,$ and using the fact that $ \xi_{i}, \ i=1,2,...n$ are orthonormal, it is clear that the problem (\ref{w2}) corresponds to the following linear system of ordinary differential equations (ODEs):
\begin{eqnarray*}	
	\left \{ \begin{array}{lc}
		M^{\textrm{T}} D_{n}^{\prime \prime}(t)+[N^{\textrm{T}}+P^{\textrm{\textrm{T}}}]D_{n}^{\prime}(t)+Q^{\textrm{T}} D_{n}(t) =G_{n}(t),\; \;  \mbox{for} \ \ t\in [0,T],  \\ [2pt]
		D_{n}(0)=U_n, \; \; D_{n}^{\prime}(0) = V_n,
	\end{array} \right.
\end{eqnarray*}
where  $D_{n}(t)=(d_{1,{n}}(t),d_{2,{n}}(t),...,d_{n,{n}}(t))^{\textrm T},$ the entries of the matrix $M,N,P,Q$  are
\begin{eqnarray*}
	M=\left[(\rho \ \xi_{i},\xi_{j})\right]_{n\times n}, \ \ \ N=\left[(\mu \ \xi_{i},\xi_{j})\right]_{n\times n},  \ \ \
	P=\left[(\kappa \ \xi_{i,xx},\xi_{j, xx})\right]_{n \times n},\ \
	Q=\left[(r\ \xi_{i,xx},\xi_{j,xx})\right]_{n \times n},
\end{eqnarray*}
and $G_{j}(t)= g(t)\xi_{j}(\ell), \ G_{n}(t)= (G_{1}(t),G_{2}(t),...,G_{n}(t))^{\textrm{T}},$ $U_j=(u_{0},\xi_{j}),$  $V_j=(v_0,\xi_{j}), \\ U_{n}=(U_{1,},U_{2},...,U_{n})^{\textrm{T}},$   $V_{n}=(V_{1},V_{2},...,V_{n})^{\textrm{T}}.$  \par
By the Carath\'eodory theorem for ODEs (see \cite{conddington}, Chapter 2, Theorem 1.1), for every $ n\geq 1$ there exists a unique solution $ u_{n}\in C^1([0,T];W_n)$ with $ u_{n}^{\prime \prime}\in L^{2}(0,T;W_n)$ of problem (\ref{w2}).
\begin{theorem} \label{342}
	Let assumptions ({\ref{3-5}}) hold true. Then, in the perspective of Definition \ref{d1}, there exists a unique weak solution $ u$ to the direct problem (\ref{1-0}). Moreover,
	\begin{eqnarray}
		\Vert u_{xx} \Vert^2_{L^{\infty}(0,T;L^2(0,\ell))} &\leq& \frac{2(C_{0}^2+1)}{r_0} \left[ \frac{2(1+T) \ell^3}{3r_0}\, \Vert g' \Vert^2_{L^2(0,T)}+R_{0}(u_{0},v_{0}) \right], \label{00} \\
		%\Vert u_{xx} \Vert^2_{L^2(0,T;L^2(0,\ell))} &\leq& C_1^2 \ (\exp(T)-1) \ \Vert g^\prime \Vert^2_{L^2(0,T)} \label{1.10} \\ [2pt]
		\Vert u \Vert^2_{L^2(0,T;\mathcal{V}_1^2(0,\ell))}&\leq& \frac{2 C^{\ast}C_{0}^2}{r_0} \ \left[ \frac{2(1+T) \ell^3}{3r_0}\, \Vert g' \Vert^2_{L^2(0,T)}+R_{0}(u_{0},v_{0}) \right]  ,\label{n4}\\ [2pt]
		\Vert u_{t}\Vert^2_{L^2(0,T;\mathcal{V}_1^2(0,\ell))}&\leq& \frac{ C^\ast (C_0^2+1)}{2\kappa_0} \left[ \frac{2(1+T) \ell^3}{3r_0}\, \Vert g' \Vert^2_{L^2(0,T)}+R_{0}(u_{0},v_{0}) \right] ,~~~~	\label{1.12}
	\end{eqnarray}
	and
	\begin{eqnarray}
		\Vert u_{tt} \Vert^2_{L^2(0,T;L^2(0,\ell))}&\leq& \frac{C_1^2}{2\rho_0} \left[ (1+T) \ell^3\Vert g^{\prime}\Vert^2_{L^2(0,T)}+R_{1}(u_0,v_0)\right] \label{09},
	\end{eqnarray}
	where $ R_{0}(u_0,v_0):=\rho_1 \Vert v_{0}\Vert^2_{L^2(0,\ell)}  +r_1 \Vert u_{0,xx}\Vert^2_{L^2(0,\ell)},$  $R_{1}(u_0,v_0):= \Vert v_{0}\Vert^2_{L^2(0,\ell)}  + \Vert u_{0,xx}\Vert^2_{L^2(0,\ell)}+ \Vert v_{0,xx}\Vert^2_{L^2(0,\ell)},$ $ C_0^2= \left(\exp(T)-1\right),$    $r_0, \rho_0,  \kappa_0 $ are the constants given in (\ref{3-5}), $C^{\ast}$ is from (\ref{4.8}), and the constant $ C_1>0$ is introduced in the proof.
\end{theorem}
{\bf Proof.}
Consider the Galerkin approximation of (\ref{1-1}), multiply it by  $ 2  d_{i,n}^{\prime}(t)$ and sum over $ i=1,2,3...,n.$ Further, instead of doing integration by parts as in (\ref{w2}), we use the formal identities
\begin{eqnarray}
	\left. \begin{array}{ll}
		2\big(r(x)u_{n,xx}\big)_{xx} u_{n}^{\prime} \equiv 2 \big[(r(x)u_{n,xx})_x u_n^{\prime}-r(x)u_{n,xx} u_{n,x}^{\prime}\big]_x+\left (r(x)u_{n,xx}^2\right )^{\prime}, \label{i1} \\ [2.5pt]
		2\big(\kappa(x)u_{n,xx}^{\prime}\big)_{xx} u_n^{\prime} \equiv 2 \big[(\kappa(x)u_{n,xx}^{\prime})_x u_n^{\prime}-\kappa(x)u_{n,xx}^{\prime} u_{n,x}^{\prime}\big]_x+2\kappa(x)(u_{n,xx}^{\prime})^2\label{i2}.
	\end{array} \right.
\end{eqnarray}
Integrating by parts, using the initial and boundary conditions of (\ref{1-1}), we obtain the following \emph{energy identity}:
\begin{eqnarray}\label{1-9}
	\lefteqn{\int_0^\ell \Big (\rho(x) u_{n}^{\prime }(t)^2+r(x)u_{n,xx}^2(t)\Big)dx +2\int_0^t \int_0^\ell \mu (x) (u_{n}^{\prime})^2 dx d\tau  +2\int_0^t \int_0^\ell \kappa(x)(u_{n,xx}^{\prime})^2 dx  d \tau}   \nonumber \\
	&&= 2g(t)u_{n} (\ell,t)- 2 \int_0^t g'(\tau)u _{n}(\ell,\tau) +  \int_{0}^{\ell}\Big(\rho(x)v_{0,n}^2+r(x)(u_{0,n})_{xx}^2\Big)dx.~~~
\end{eqnarray}
We employ the $\varepsilon$-inequality $2ab \leq (1/\varepsilon)\, a^2+\varepsilon \,b^2, \varepsilon>0$ frequently in the proof. Apply this inequality in the first two terms of the right-hand side of (\ref{1-9}), and then use the trace inequalities
\begin{eqnarray}
	u_{n}^2(\ell,t) &\leq& \frac {\ell^3}{3}  \int_0^\ell u_{n,xx}^2(x,t) dx
	,   \label{1-10} \\[2.5pt]
	g^2(t) &\leq& T \Vert g^{\prime}\Vert^2_{L^2(0,T)}, \ \ \ \mbox{for all}\ t \in [0,T],\label{12}
\end{eqnarray}
where the inequality (\ref{1-10}) is a consequence of the identity
\begin{eqnarray*}
	u_{n}(\ell,t) \equiv \int_0^\ell (\ell-x) u_{n,xx}(x,t) dx, \ \mbox{for all}\ t \in [0,T].
\end{eqnarray*}
Then we choose the arbitrary constant $\varepsilon>0$ from the condition
$r_0-\ell^3\,\varepsilon/3>0$ as follows: $\varepsilon=3r_0/(2\ell^3)$. After elementary transformations, we obtain the following main integral inequality
\begin{eqnarray}\label{1-11}
	\lefteqn{	\rho_0\int_0^\ell  u_n^{\prime}(t)^2dx + \frac{r_0}{2}\,\int_0^\ell u_{n,xx}^2(t) dx+2\int_{0}^t\int_0^\ell \mu(x) (u_{n}^{\prime})^2 dx  d \tau   +2\int_{0}^t\int_0^\ell  \kappa(x)(u_{n,xx}^{\prime})^2 dx  d \tau} \nonumber \\
	&& \le \frac{r_0}{2}\, \int_0^t \int_0^\ell u_{n,xx}^2dx  d \tau 	 + \frac{2\ell^3 (1+T)}{3r_0} \, \Vert g' \Vert^2_{L^2(0,T)} +\rho_1 \Vert v_{0,n}\Vert^2_{L^2(0,\ell)} +r_1 \Vert u_{0,nxx}\Vert^2_{L^2(0,\ell)}.
\end{eqnarray}
The first consequence of the integral inequality (\ref{1-11}) is that
\begin{eqnarray*}
	\int_0^\ell u_{n,xx}^2(t) \ dx &\le& \int_0^t \int_0^\ell u_{n,xx}^2 \ dx d \tau+ \frac{2}{r_0} R(g,u_0,v_0),
\end{eqnarray*}
where $$ R(g,u_{0},v_{0}):=\frac{2(1+T) \ell^3}{3r_0}\, \Vert g' \Vert^2_{L^2(0,T)} +\rho_1 \Vert v_{0}\Vert^2_{L^2(0,\ell)}+r_1 \Vert u_{0,xx}\Vert^2_{L^2(0,\ell)}.$$
By invoking the Gr\"onwall-Bellmann inequality, we obtain
\begin{eqnarray}\label{1-12}
	\int_0^\ell u_{n,xx}^2 (t)\ dx &\leq& \frac{2}{r_0}\exp(t) R(g,u_0,v_0).
\end{eqnarray}
Integrating inequality (\ref{1-12}) over $[0,T],$ we arrive at the first required estimate as follows:
\begin{eqnarray}\label{1.2}
	\Vert u_{n,xx}\Vert^2_{L^2(0,T;L^2(0,\ell))}&\leq& \frac{2C_{0}^2}{r_0} R(g,u_0,v_0),
\end{eqnarray}
where $ C_{0}^2=\left(\exp(T)-1\right).$ Taking maximum over $ t\in[0,T]$ in (\ref{1-12}), we get
\begin{eqnarray}\label{0}
	\max_{t\in[0,T]}\Vert u_{n,xx}(t)\Vert^2_{L^2(0,\ell)}&\leq&\frac{2( C_0^2+1)}{r_0} R(g,u_0,v_0) .
\end{eqnarray}
Since $ \Vert u_{n}(t)\Vert^2_{\mathcal{V}_{1}^2(0,\ell)}\leq C^{\ast} \Vert u_{n,xx}(t)\Vert^2_{L^2(0,\ell)}$ by (\ref{4.8}), using the estimate (\ref{1.2}), we obtain
\begin{eqnarray}\label{1.3}
	\Vert u_{n}\Vert^2_{L^2(0,T;\mathcal{V}_{1}^2(0,\ell))}&\leq& \frac{2 C^{\ast}C_{0}^2}{r_0} R(g,u_0,v_0).
\end{eqnarray}
The second consequence of (\ref{1-11}) and (\ref{1.2}) is the inequality
\begin{eqnarray}\label{2.98}
	2\int_0^t \int_0^\ell\kappa(x) (u_{n,xx}^\prime)^2 dx d \tau &\le&  \frac{r_0}{2}\, \int_0^t \int_0^\ell u_{n,xx}^2dx d \tau+ R(g,u_0,v_0) \\ \nonumber 
	&\leq &(C_{0}^2+1) R(g,u_0,v_0) .
\end{eqnarray}
The estimate (\ref{2.98}) and again the equality of norms lead to the inequality
\begin{eqnarray}\label{1.4}
	\Vert u_{n}^{\prime} \Vert^2_{L^2(0,T;\mathcal{V}_1^2(0,\ell))}&\leq& \frac{C^\ast(C_{0}^2+1)}{2\kappa_0} R(g,u_0,v_0).
\end{eqnarray}

	To determine the estimate $ \Vert u_{n}^{\prime \prime}\Vert^2_{L^2(0,T;L^2(0,\ell))},$ we proceed as follows.
Multiply the Galerkin approximation of (\ref{1-1}) by $2 d_{i,n}^{\prime \prime}(t)$ and use the following formal identities
\begin{eqnarray*}
	2\big(r(x)u_{n,xx}\big)_{xx} u_{n}^{\prime \prime} &\equiv& 2 \big[(r(x)u_{n,xx})_x u_{n}^{\prime \prime}-r(x)u_{n,xx} u_{n,x}^{\prime \prime}\big]_x +2\left (r(x)u_{n,xx} u_{n,xx}^{\prime \prime}\right),\\
	2\big(\kappa(x)u_{n,xx}^{\prime}\big)_{xx} u_{n}^{\prime \prime}& \equiv& 2 \big[(\kappa(x)u_{n,xx}^{\prime})_x u_{n}^{\prime \prime}-\kappa(x)u_{n,xx}^{\prime} u_{n,x}^{\prime \prime}\big]_x+\left (\kappa(x)(u_{n,xx}^{\prime})^2\right )^{\prime}.
\end{eqnarray*}
Integrate by parts and invoking the initial and boundary conditions of (\ref{1-1}), we obtain the following second energy identity
\begin{eqnarray}\label{1-14}
	\lefteqn{\int_{0}^{\ell} \mu(x)u_{n}^{\prime}(t)^2 dx+ \int_0^\ell \kappa(x)\left(u_{n,xx}^\prime(t)\right)^2 dx+ 2 \int_{0}^{t}\int_{0}^{\ell} \rho(x)\left(u_{n}^{\prime \prime}\right)^2 dx    d\tau} \nonumber \\
	&=& 2\int_{0}^{t} \int_{0}^{\ell}r(x)(u_{n,xx}^{\prime})^2 dx  d\tau- 2 \int_0^{\ell}r(x)u_{n,xx}(t)u_{n,xx}^{\prime}(t)dx +2g(t)u_{n}^{\prime}(\ell,t) - 2\int_{0}^{t} g^{\prime}(\tau)u_{n}^{\prime}(\ell,\tau) d\tau  \nonumber \\ &&
	+ \int_{0}^{\ell}\kappa(x)v_{0,nxx}^2dx  +\int_{0}^{\ell}\mu(x)v_{0,n}^2dx +2\int_{0}^{\ell} r(x)u_{0,nxx}v_{0,nxx}dx.
\end{eqnarray}
By applying Cauchy's inequality to the second, third, fourth and seventh terms on the right-hand side of (\ref{1-14}), using the inequality (\ref{12}) and
\begin{eqnarray}\label{mnbv}
	(u_{n}^{\prime}(\ell,t))^2&\leq& \frac{\ell^3}{3}\int_{0}^{\ell}(u_{n,xx}^{\prime})^2 dx,
\end{eqnarray}
we get,
\begin{eqnarray*}
	\lefteqn{\int_{0}^{\ell} \mu(x)u_{n}^{\prime}(t)^2 dx+ \int_0^\ell \kappa(x)\left(u_{n,xx}^\prime(t)\right)^2 dx+ 2 \int_{0}^{t}\int_{0}^{\ell} \rho(x)(u_{n}^{\prime \prime})^2 dx d\tau} \nonumber \\
	&\leq&\left( 2r_1+\frac{\ell^3\epsilon}{3}\right)\int_{0}^{t} \int_{0}^{\ell}(u_{n,xx}^{\prime})^2 dx  d\tau+\frac{r_1^2}{\epsilon} \int_{0}^{\ell} u_{n,xx}^2(t)dx  +\left(\frac{\ell^3\epsilon}{3}+\epsilon\right)\int_{0}^{\ell} (u_{n,xx}^{\prime}(t))^2dx \ \nonumber \\&&+\frac{1+T}{\epsilon}\int_{0}^t g^{\prime}(\tau)^2d\tau +\mu_1 \Vert v_{0}\Vert^2_{L^2(0,\ell)} \nonumber  +\left(\kappa_1+1\right)\Vert v_{0,xx}\Vert^2_{L^2(0,\ell)}  +r_1^2\Vert u_{0,xx}\Vert^2_{L^2(0,\ell)}.
\end{eqnarray*}
Now choose $ \epsilon =3\kappa_0/2(\ell^3+3)$ from the condition $ \kappa_0-\left(\ell^3 \epsilon/3+\epsilon\right)>0$ and invoking the estimates (\ref{0}), (\ref{2.98}), we obtain
\begin{eqnarray*}\label{1-15}
	2 \int_{0}^{t}\int_{0}^{\ell} \rho(x)\left(u_{n}^{\prime \prime}\right)^2 dx d\tau &\leq& R_2\ell^3(1+T)\Vert g^\prime\Vert^2_{L^2(0,T)}+ R_3 \Vert v_{0}\Vert^2_{L^2(0,\ell)}\nonumber \\ &&+R_4\Vert u_{0,xx}\Vert^2_{L^2(0,\ell)}+\left(\kappa_1+1\right)\Vert v_{0,xx}\Vert^2_{L^2(0,\ell)},~~~~
\end{eqnarray*}	
where \begin{eqnarray*}
	R_1&=&\left(\frac{4r_1^2(\ell^3+3)}{3 r_0}+\frac{\kappa_0\ell^3}{4(\ell^3+3)}+r_1\right)\frac{(C_0^2+1)}{\kappa_0}, \ \
	R_2=R_1\frac{2}{3 r_0}+\frac{2(l^3+3)}{3l^3\kappa_0},  \\
	R_3&=&R_1\rho_1+\mu_1, \ \ \ R_4=R_1r_1+r_1^2.
\end{eqnarray*}
Choosing $ C_1^2=\max\{R_2,R_3,R_4,\kappa_1+1\}, $ we get
\begin{eqnarray}\label{1.7}
	\Vert u_{n}^{\prime \prime}\Vert_{L^2(0,T,L^2(0,\ell))}&\leq& \frac{C_1^2}{2\rho_0} \Big[\ell^3(1+T)\Vert g^\prime\Vert^2_{L^2(0,T)} + \Vert v_{0}\Vert^2_{L^2(0,\ell)}+\Vert u_{0,xx}\Vert^2_{L^2(0,\ell)} +\Vert v_{0,xx}\Vert^2_{L^2(0,\ell)}\Big].
\end{eqnarray}
Consequently, by using the estimates (\ref{1.3}), (\ref{1.2}), (\ref{1.4}) and (\ref{1.7}), we obtain that the sequences $\{u_n\},\ \{u_{n,xx}\},$  \\ $\{u_{n,xx}^\prime\}, \{u_n^\prime\},$ and $ \{u_{n}^{\prime \prime}\}$ are bounded in $ L^{2}(0,T;\mathcal{V}_1^2(0,\ell)),$ $ L^{2}(0,T;L^2(0,\ell)),$ $ L^2(0,T;L^2(0,\ell)),$ \\ $ L^2(0,T;\mathcal{V}_1^2(0,\ell)),$  $L^2(0,T;L^2(0,\ell))$ respectively. We now use the Banach-Alaoglu weak compactness theorem (see,  Theorem 3.16, \cite{brezis2010functional}) to deduce that there exists a subsequence $ \{u_{n_k}\}$ of $ u_{n}$ and functions $ u\in L^2(0,T;\mathcal{V}_1^2(0,\ell)),$ $ u_{xx}\in L^2(0,T;L^2(0,\ell)),$ $ u_{xx}^{\prime}\in L^2(0,T;L^2(0,\ell)),$ $ u^{\prime} \in L^2(0,T;\mathcal{V}_1^2(0,\ell)),$  and $u^{\prime \prime} \in L^2(0,T;L^2(0,\ell))$ such that
\begin{eqnarray}\label{p}
	\left \{ \begin{array}{lclccl}
		u_{n_{k}}&\rightharpoonup& u \ &\mbox{weakly  in}  &  L^{2}(0,T;\mathcal{V}_1^{2}(0,\ell)) \\
		u_{n_k,xx}&\rightharpoonup& u_{xx}\ &\mbox{weakly  in} & L^{2}(0,T;L^{2}(0,\ell))\\
		u_{n_k,xx}^{\prime}&\rightharpoonup& u_{xx}^{\prime}\ &\mbox{weakly  in} & L^{2}(0,T;L^{2}(0,\ell))\\
		u_{n_{k}}^{\prime}&\rightharpoonup& u^{\prime} \ &\mbox{weakly in} & L^{2}(0,T;\mathcal{V}_1^{2}(0,\ell))\\
		u_{n_{k}}^{\prime \prime} &\rightharpoonup& u^{\prime \prime} \ &\mbox{weakly   in} & L^{2}(0,T;L^{2}(0,\ell)).
	\end{array} \right.
\end{eqnarray}
We should be able to derive the weak solution $ u $ of the direct problem (\ref{1-1}) by passing the limit on the weak form (\ref{w2}). The solution $ u$ also satisfies estimates (\ref{00})-(\ref{09}). \par
The uniqueness of weak solution of direct problem (\ref{1-1}) can be proved by using (\ref{n4}) and (\ref{1.12}). Suppose that there are two weak solutions $ u_1$ and $u_2$ in $ \mathcal{V}_1^2(0,\ell)$ of the direct problem (\ref{1-1}). Then a function $ \mathcal{U}(x,t)=u_1(x,t)-u_2(x,t)$ that solves the direct problem (\ref{1-1}) with homogeneous initial and boundary conditions. The estimates (\ref{n4}) and (\ref{1.12}) applied to this problem imply that $ \Vert \mathcal{U} \Vert_{L^{\infty}(0,T;\mathcal{V}_1^2(0,\ell))}=0.$ Hence the homogeneity of initial and boundary conditions imply that $ \mathcal{U}(x,t)=0$, $\forall \ (x,t)\in (0,\ell)\times(0,T).$\par
It is still necessary to verify that $ u(t)$ satisfy the initial condition $ u(0)=u_0$ and $ u_{t}(0)=v_0.$ Taking $ u \in C([0,T];\mathcal{V}_1^2(0,\ell))$ and $ u_ {t} \in C([0,T];L^2(0,\ell))$ into account, choosing a test function $ v\in C^2([0,T];\mathcal{V}_1^2(0,\ell))$ with $ v(T)=0$ and $ v^\prime(T)=0$  and arguing as in \cite{anjuna:2021} or Theorem 3 in \cite{Baysal:Hasanov:2019}, one can verify the initial data. This completes the proof. \hfill$\Box$
\subsection{Regularity of Weak Solutions}
In this subsection, we study the regularity of the weak solution, which is required to show the compactness of the input output operator and the Fr\'echet derivative of the functionals.  For simplicity, we take $u_0=v_0\equiv 0.$
\begin{theorem}\label{8}
	Let conditions (\ref{3-5}) hold. Assume that  the following regularity and consistency conditions are also satisfied:
	\begin{eqnarray} \label{3-2}
		\left \{ \begin{array}{ll}
			r,\ \kappa\in H^2(0,\ell), \, \Vert r \Vert_{H^2(0,\ell)}\leq r_2, \, \Vert \kappa \Vert_{H^2(0,\ell)}\leq \kappa_2 \nonumber \\
			g\in H^2(0,T),\ g(0)=0,\ g^\prime(0)=0.
		\end{array}\right.
	\end{eqnarray}
	Then the following estimate holds for regular weak solution of (\ref{1-1}) with enhanced regularity $ u\in H^1(0,T;H^4(0,\ell)),$ $u_{tt}\in L^2(0,T;\mathcal{V}_1^2(0,\ell)),$  $ u_{ttt}\in L^2(0,T;L^2(0,\ell)):$
	\begin{eqnarray}\label{4.7}
		\Vert u_{ttt}\Vert^2_{L^2(0,T;L^2(0,\ell))}\leq \frac{C_{5}^2}{2\rho_0}\exp(C_5^2T)\Vert g\Vert^2_{H^2(0,T)},
	\end{eqnarray}
	where
	\begin{eqnarray*}
		C_5^2:=\frac{2}{\kappa_0}\max\left(\frac{2(1+T)(\ell^3+3)}{3\kappa_0}, \, \left[\frac{\kappa_0\ell^3}{2(\ell^3+3)}+2r_1+\frac{2Tr_1^2 (\ell^3+3)}{3\kappa_0}\right]\right ).
	\end{eqnarray*}
\end{theorem}
{\bf Proof.}
Multiply equation (\ref{1-1}) by $2 u_{xxxxt},$ integrate over $ (0,\ell)\times(0,t)$ and apply Cauchy's $ \epsilon$-inequality with $ \epsilon=\kappa_0/6,$ where $ \kappa_0$ is given in (\ref{3-5}),  we arrive at
\begin{eqnarray}\label{4.15}
	\lefteqn{\int_0^\ell r(x)u_{xxxx}^2(t)dx+\kappa_0\int_0^t\int_0^{\ell} u_{xxxx\tau}^2 dx d\tau} \nonumber \\&\leq& \frac{6}{\kappa_0}\Big[ \int_0^t \int_0^\ell \rho^2u_{\tau \tau}^2 dx d\tau+ \int_0^t \int_0^\ell \mu^2 u_{ \tau}^2 dx d\tau
	+4 \int_0^t \int_0^\ell (r^{\prime})^2u_{xxx}^2 dx d\tau  +4 \int_0^t \int_0^\ell (\kappa^{\prime})^2 u_{xxx\tau}^2 dx d\tau\nonumber \\&&+\int_0^t\int_0^\ell (r^{\prime \prime})^2u_{xx}^2dx d\tau +\int_0^t\int_0^\ell (\kappa^{\prime \prime})^2 u_{xx\tau}^2 dx d\tau\Big]  =:\,\frac{6}{\kappa_0} \sum_{i=1}^{6} I_{i}.
\end{eqnarray}
Notice that,
\begin{eqnarray}\label{5.78}
	I_3\leq 4 \int_{0}^{t}\!\!\!\Vert r^{\prime}\Vert^2_{L^{\infty}(0,\ell)}\Vert u_{xxx}(\tau)\Vert^2_{L^2(0,\ell)} d\tau 
	\leq4 \bar{C} \Vert r^{\prime}\Vert^2_{H^1(0,\ell)} \Vert u_{xxx}\Vert^2_{L^2(0,T;L^2(0,\ell))},~~
\end{eqnarray}
and similarly, $ I_{4}\leq 4 \bar{C} \kappa_{2}^2 \Vert u_{xxxt}\Vert^2_{L^2(0,T;L^2(0,\ell))},$
where we used the fact that $ H^1(0,\ell)$ is continuously embedded in $ L^\infty(0,\ell),$ that is, $$ \Vert u \Vert_{L^{\infty}(0,\ell)}\leq \bar{C}(\ell)\Vert u\Vert_{H^1(0,\ell)}, \ \  \bar{C}(\ell)=\sqrt{2}\max\{ 1/\sqrt{\ell},\ell\},$$ (see, \cite{salsa2015}, Section 7.10.4). Further,  the integrals $I_5$ and $I_6$ can be estimated as follows
	\begin{eqnarray}
	I_5&\leq& \int_0^t \Vert u_{xx}(\tau)\Vert^2_{L^\infty(0,\ell)} \Vert r^{\prime \prime}\Vert^2_{L^2(0,\ell)}d\tau  \leq\bar{C}  r_2^2 \Vert u \Vert^2_{L^2(0,t;H^3(0,\ell))}\nonumber \\&=& \bar{C} r_2^2\bigg(\Vert u \Vert^2_{L^2(0,t;\mathcal{V}^2_1(0,\ell))} \! + \! \Vert u_{xxx}\Vert^2_{L^2(0,t;L^2(0,\ell))}\bigg),\label{3.9}\\
	I_6&\leq& \bar{C} \kappa_2^2 \Vert u_{t}\Vert^2_{L^2(0,t;H^3(0,\ell))}= \bar{C}\kappa_2^2\bigg(\Vert u_t \Vert^2_{L^2(0,t;\mathcal{V}_1^2(0,\ell))} \!+ \!\Vert u_{xxx\tau} \Vert^2_{L^2(0,t;L^2(0,\ell))}\bigg). \ \ \label{3.10} \ \
\end{eqnarray}
From the equality of norms (see,  \cite{brezis2010functional}, page 217), for every integer $ j,$ $ 1\leq j\leq m-1$ and for every $ \lambda >0,$ there exists a constant $ \tilde{C}(\lambda,\ell)$ such that $$ \Vert D^ju\Vert_{L^p(0,\ell)}\leq \lambda \Vert D^m u \Vert_{L^{p}(0,\ell)} + \tilde{C}(\lambda,\ell) \Vert u \Vert_{L^{p}(0,\ell)}, \ \forall u \in W^{m,p}(0,\ell),\ \ p\geq 1.$$
By invoking this result, let us estimate the integrals $\Vert u_{xxx}\Vert^2_{L^2(0,T;L^2(0,\ell))}$ and $\Vert u_{xxxt}\Vert^2_{L^2(0,T;L^2(0,\ell))}.$  For any $ \lambda_1>0,$ we obtain that
\begin{eqnarray}\label{ght}
	\int_0^t\int_0^\ell u_{xxx\tau}^2 dx d\tau &\leq& 2\lambda_1 \int_0^t\int_0^\ell u_{xxxx\tau}^2 dx d\tau +\tilde{C}(\lambda_1,\ell)\int_0^t\int_0^\ell u_{\tau}^2 dx d\tau,\label{3.7}
\end{eqnarray}
and similar estimate holds for $ \Vert u_{xxx}\Vert^2_{L^2(0,T;L^2(0,\ell))}$ with $ \lambda_2>0.$
Using the above estimates (\ref{5.78})-(\ref{ght}) in (\ref{4.15}) and choosing  $ \lambda_1=\kappa_0^2/(120\kappa_2^2\bar{C})$, $\lambda_2=1/2,$ we obtain
\begin{eqnarray}
	\sum_{i=1}^{6}I_{i}&\leq& \frac{6}{\kappa_0} \left(\rho_1^2 \Vert u_{tt}\Vert^2_{L^2(0,T;L^2(0,\ell))} +\left(\mu_1^2+\bar{C} \kappa_2^2(1+5\tilde{C})\right)\Vert u_{t}\Vert^2_{L^2(0,T;\mathcal{V}_1^2(0,\ell))} \right. \nonumber \\ && \left.+\bar{C}r_2^2(1+5\tilde{C}) \ \Vert u \Vert^2_{L^2(0,T;\mathcal{V}_{1}^2(0,\ell))} +  5\bar{C} r_2^2\int_0^{t}\int_0^{\ell} u_{xxxx}^2 \ dxd\tau \right) \nonumber \\ &&+\frac{\kappa_0}{2}\int_0^t \int_0^\ell u_{xxxx\tau}^2 \ dxd\tau.
\end{eqnarray}
Making use of the estimates (\ref{n4})-(\ref{09}), we further have
\begin{eqnarray}\label{3.18}
	r_0\int_0^\ell u_{xxxx}^2(t)dx&+&\frac{\kappa_0}{2}\int_0^t\int_0^{\ell}u_{xxxx\tau}^2 dx d\tau   \leq  C_2^2 \Vert g^{\prime}\Vert^2_{L^2(0,T)}+\frac{30\bar{C}r_2^2}{\kappa_0}\int_0^t\int_0^\ell u_{xxxx}^2 dx d\tau,
\end{eqnarray}
with the constant
\begin{eqnarray*} C_2^2=\frac{6 \ell^3(1+T)}{\kappa_0}  \times\max \left ( \frac{\rho_1^2C_1^2}{2\rho_0},\frac{ 4\bar{C}r_2^2(1+5\tilde{C}) C^{\ast}C_{0}^2}{3r^2_0}, \left [ \mu_1^2+(1+5\tilde{C})\bar{C}\kappa_2^2\right ]	\frac{C^\ast (C_0^2+1)}{3r_0\kappa_0}\right ).
\end{eqnarray*}
Applying Gr\"onwall's inequality and then integrating over $(0,T)$ , we get
\begin{eqnarray}\label{3.15}
	\Vert u_{xxxx} \Vert^2_{L^2(0,T;L^2(0,\ell))}\leq C_3^2 \Vert g^\prime\Vert^2_{L^2(0,T)},
\end{eqnarray}
where $C_3^2=\frac{C_2^2\kappa_0}{30\bar{C}r_2^2} \left[\exp\Big(\frac{30\bar{C}r_2^2}{\kappa_0r_0}T\Big)-1\right].$
Next, let us note that \begin{eqnarray}
	\Vert u \Vert^2_{L^2(0,T;H^4(0,\ell))}&=& \Vert u \Vert^2_{L^2(0,T;H^3(0,\ell))}+\Vert u_{xxxx}\Vert^2_{L^2(0,T;L^2(0,\ell))},\label{3.12} \\
	\Vert u_{t}\Vert^2_{L^2(0,T;H^4(0,\ell))}&=& \Vert u_t \Vert^2_{L^2(0,T;H^3(0,\ell))}+\Vert u_{xxxxt} \Vert^2_{L^2(0,T;L^2(0,\ell))}\label{3.13}.
\end{eqnarray}
To estimate the first term on the right-hand side of (\ref{3.12}), we use Ehrling's lemma (see, \cite{MRenardyRogers}, Theorem 7.30), that for any $ \lambda_3>0, $ there exist $ C(\lambda_3)$ such that for any $  u \in {L^2(0,T;H^4(0,\ell))},$ we get
\begin{eqnarray*}\label{3.17}
	\Vert u \Vert^2_{L^2(0,T;H^3(0,\ell))}\leq \lambda_3 \Vert u \Vert^2_{L^2(0,T;H^4(0,\ell))}+C(\lambda_3) \Vert u \Vert^2_{L^2(0,T;L^2(0,\ell))}.
\end{eqnarray*}
Choosing $ \lambda_3=\frac{1}{2},$ substituting this into (\ref{3.12}) and using (\ref{3.15}), \ref{n4}, we obtain
\begin{eqnarray}\label{poi}
	\Vert u \Vert^2_{L^2(0,T;H^4(0,\ell))}\leq 2\left( C_3^2+\frac{4CC^\ast C_0^2}{3r_0^2}(1+T)\ell^3\right) \Vert g^\prime \Vert^2_{L^2(0,T)}.
\end{eqnarray}
The second consequence of (\ref{3.18}) and (\ref{3.15}) is the following
\begin{eqnarray}
	\Vert u_{xxxxt}\Vert^2_{L^2(0,T;L^2(0,\ell))}\leq\frac{2}{\kappa_0}\left(C_2^2+ \frac{30\bar{C}C_3^2 r_2^2}{\kappa_0}\right) \Vert g^\prime\Vert^2_{L^2(0,T)}.
\end{eqnarray}
Again invoking Ehrling's lemma for $ u_{t}\in L^2(0,T;H^4(0,\ell)),$  (\ref{3.13}) and \ref{1.12}, we get
\begin{eqnarray}\label{pou}
	\Vert u_t \Vert^2_{L^2(0,T;H^4(0,\ell))}\leq C_4^2 \Vert g^\prime \Vert^2_{L^2(0,T)},
\end{eqnarray}
where $C_4^2=\frac{2}{\kappa_0}\left(2\left(C_2^2+ \frac{30\bar{C}C_3^2 r_2^2}{\kappa_0}\right)+\frac{C C^\ast (C_0^2+1)\ell^3(1+T)}{3r_0}\right).$ Consequently, from (\ref{poi}) and (\ref{pou}), we infer that $ u \in H^1(0,T;H^4(0,\ell)).$

To estimate $ \Vert u_{ttt}\Vert^2_{L^2(0,T;L^2(0,\ell))}$ and $ \Vert u_{xxtt}\Vert^2_{L^2(0,T;L^2(0,\ell))},$ we proceed as follows.
Formally differentiate  (\ref{1-1}) with respect to time and multiply by $ 2 u_{ttt},$ use the crucial identities
\begin{eqnarray*}\label{9}
	\left. \begin{array}{ll}
		2\big(r(x)u_{xx}\big)_{xxt} u_{ttt} \equiv 2 \big[(r(x)u_{xx})_{xt} u_{ttt}-r(x)u_{xxt} u_{xttt}\big]_x+2r(x)u_{xxt} u_{xxttt}, \nonumber \\
		2\big(\kappa(x)u_{xxt}\big)_{xxt} u_{ttt} \equiv 2 \big[(\kappa(x)u_{xxt})_{xt} u_{ttt}-\kappa(x)u_{xxtt} u_{xttt}\big]_x+\kappa(x)(u_{xxtt}^{2})_t,
	\end{array} \right.
\end{eqnarray*}
and integrate over $(0,\ell)\times(0,t), t\in(0,T).$ Then integrating by parts using the initial and boundary conditions of (\ref{1-1}), we obtain
\begin{eqnarray}\label{1-r}
	\lefteqn{	\int_0^\ell \Big (\mu(x) u_{tt}(t)^2+\kappa(x)u_{xxtt}(t)^2\Big)dx +2\int_0^t \int_0^\ell \rho (x) u_{\tau \tau \tau}^{2} dx d\tau}  \nonumber \\ &=&-2\int_0^t \int_0^\ell r(x)(u_{xx\tau} u_{xx\tau \tau \tau})\ dx  d \tau +2\int_0^t g^{\prime}(\tau)u_{\tau \tau \tau } (\ell,\tau) \ d\tau  +\int_0^\ell \mu(x)u_{tt}^2(x,0^{+})dx \nonumber \\
	&&+\int_0^\ell \kappa(x)u_{xxtt}^2(x,0^{+}) dx.
\end{eqnarray}
Let us evaluate the third and fourth right-hand side integrals by using the initial data given in (\ref{1-1}). First, we deduce that
\begin{eqnarray*}
	\int_0^\ell \mu(x)u_{tt}^2(x,0^+)dx&=&\int_0^\ell \frac{\mu(x)}{\rho^2(x)}\Big(\mu(x)u_{t}(x,0^{+})+(r(x)u_{xx}(x,0^{+}))_{xx} +(\kappa(x)u_{xxt}(x,0^+))_{xx}\Big)^2 dx=0,
\end{eqnarray*}
since $u_{t}(x,0^+)=u_{xx}(x,0^+)=u_{xxt}(x,0^+)=0.$
Similarly, we also obtain that
$\int_0^\ell \kappa(x)u_{xxtt}^2(x,0^+)dx=0.$
Now the integration by parts with respect to time over the first two right-hand side integrals of (\ref{1-r}) leads to the following
\begin{eqnarray}\label{1-e}
	\lefteqn{-2\int_0^t \int_0^\ell r(x)(u_{xx\tau} u_{xx\tau \tau \tau})\ dx  d \tau +2\int_0^t g^{\prime}(\tau)u_{\tau \tau \tau } (\ell,\tau) \ d\tau}\nonumber \\  &=&2 g^{\prime}(t)u_{tt}(\ell,t)-2\int_0^t g^{\prime \prime}(\tau)u_{\tau \tau}(\ell,\tau)d\tau+2\int_0^t \int_0^\ell r(x)u_{xx\tau\tau}^2dx d\tau-2\int_0^\ell r(x)u_{xxt}(t)u_{xxtt}(t) dx.
\end{eqnarray}
By substituting the identity (\ref{1-e}) into (\ref{1-r}), applying the trace inequalities (\ref{1-10}), (\ref{12}) and Cauchy's inequality, we obtain that
\begin{eqnarray*}
	\lefteqn{\int_0^{\ell} \kappa(x)u_{xxtt}(t)^2dx+2\int_{0}^t\int_0^\ell \rho(x)u_{\tau \tau \tau}^2 dx d\tau  
	}
	\nonumber \\
	&\leq& \frac{1+T}{\epsilon}\int_0^t g^{\prime \prime}(\tau)^2 d\tau+\Big(\frac{\ell^3\epsilon}{3}+\epsilon\Big)\int_0^\ell u_{xxtt}(t)^2dx  +\Big(\frac{\ell^3 \epsilon }{3}+2r_1\Big)\int_0^t\int_0^\ell u_{xx\tau\tau}^2 \ dxd\tau+\frac{r_1^2}{\epsilon}\int_0^\ell u_{xxt}(t)^2 dx.~~~
\end{eqnarray*}
Taking $ \epsilon =\frac{3\kappa_0}{2(\ell^3+3)}$ and employing the trace inequality $ u_{xxt}^2(t)\leq T \int_{0}^t u_{xx\tau \tau}^2 d\tau$ give
\begin{eqnarray}\label{1.j}
	\frac{\kappa_0}{2}\int_0^{\ell} u_{xxtt}(t)^2dx+2\rho_0\int_{0}^t\int_0^\ell u_{\tau \tau \tau}^2 dx d\tau  \leq\frac{2(1+T)(\ell^3+3)}{3\kappa_0}\int_0^t g^{\prime \prime}(\tau)^2 d\tau \nonumber \\ 
	\quad\quad+\left(\frac{\kappa_0\ell^3}{2(\ell^3+3)}+2r_1+\frac{2Tr_1^2 (\ell^3+3)}{3\kappa_0}\right)\int_0^t\int_0^\ell u_{xx\tau\tau}^2 \ dxd\tau .
\end{eqnarray}
By setting $C_5^2:=\frac{2}{\kappa_0}\max\left(\frac{2(1+T)(\ell^3+3)}{3\kappa_0}, \, \left[\frac{\kappa_0\ell^3}{2(\ell^3+3)}+2r_1+\frac{2Tr_1^2 (\ell^3+3)}{3\kappa_0}\right]\right ),$ and applying Gr\"onwall's inequality, we get
$
\int_0^{\ell} u_{xxtt}(t)^2dx\leq C_5^2 \Vert g^{\prime \prime}\Vert^{2}_{L^2(0,T)}\exp(C_5^2 t).
$
This implies that
\begin{eqnarray}\label{1.t}
	\Vert u_{xxtt}\Vert^2_{L^2(0,T;L^2(0,\ell))}\leq \Vert g^{\prime \prime}\Vert^{2}_{L^2(0,T)}\left[\exp(C_5^2 T)-1\right] .
\end{eqnarray}
Substituting (\ref{1.t}) into (\ref{1.j}), we conclude the proof of the estimate (\ref{4.7}). \hfill$\Box$
\begin{theorem}\label{t1}
	Let the conditions of Theorem \ref{8} hold. Additionally assume that,  the input $ g(t)$ and the coefficients meet the following regularity and consistency conditions:
	\begin{eqnarray} \label{3-9}
		\left \{ \begin{array}{ll}
			r,\ \kappa\in H^4(0,\ell), \,\Vert r \Vert_{H^4(0,\ell)}\leq r_3, \, \Vert \kappa \Vert_{H^4(0,\ell)}\leq \kappa_3 \nonumber \\
			\rho,\ \mu \in H^2(0,\ell), \,\Vert \rho \Vert_{H^2(0,\ell)}\leq \rho_2, \, \Vert \mu \Vert_{H^2(0,\ell)}\leq \mu_2\nonumber \\
			g\in H^3(0,T),\, g(0)=g^\prime(0)= g^{\prime \prime}(0)=0.
		\end{array}\right.
	\end{eqnarray}
	Then for the regular weak solution with improved regularity, we have the enhanced regularity  $ u \in H^1(0,T;H^6(0,\ell)),$  $ u_{tt}\in L^2(0,T;H^4(0,\ell)), \  u_{ttt} \in L^2(0,T;\mathcal{V}_1^2(0,\ell)), \ u_{tttt} \in L^2(0,T;L^2(0,\ell)), $ and the estimate
	\begin{eqnarray}\label{4.1}
		\Vert u_{tttt}\Vert^2_{L^2(0,T;L^2(0,\ell))}\leq \frac{C_{5}^2}{2\rho_0}\exp(C_5^2 T)\Vert g\Vert^2_{H^3(0,T)}.
	\end{eqnarray}
\end{theorem}
{\bf Proof.}
The proof of this theorem can be completed by following the lines of arguments of Theorem \ref{8}. 
\iffalse
Differentiate equation (\ref{1-1}) two times with respect to $ x,$ multiplying it by $ 2 u_{xxxxxt}$ and applying equivalence of norm and Ehrling's lemma, we get $ u \in L^2(0,T;H^6(0,\ell)), \ u_t \in L^2(0,T;H^6(0,\ell)).$

Similarly, differentiating equation (\ref{1-1}) twice with respect to time , multiplying by $ 2 u_{xxxxtt}$ and proceeding as above, one can obtain $ u_{tt}\in L^2(0,T;H^4(0,\ell)).$

Differentiating the equation (\ref{1-1}) twice with respect to time, multiplying by $ 2 u_{tttt}$ and following the same steps done for the estimate (\ref{4.7}), we obtain $ u_{ttt}\in L^2(0,T;\mathcal{V}_1^2(0,\ell))$ and (\ref{4.1}). \fi \hfill$\Box$
\begin{remark}\label{r1}
	From Theorem \ref{342} to Theorem \ref{t1} it follows that the weak and regular weak solutions of the Euler-Bernoulli beam equation with Kelvin-Voigt damping term $(\kappa (x)u_{xxt})_{xx})$ has more enhanced regularity property than corresponding weak solutions of this equation without this term.
\end{remark}

\section{Solvability of Regularized Inverse Problems}\label{s5}

In this section, using the regularity of the solution to the direct problem (\ref{1-1}), we prove the compactness as well as the Lipschitz continuity of the input-output operators $ \Phi$ and $ \Psi$. The lower semi-continuity of the Tikhonov functionals $ J_1$ and $ J_2$ will result from the Lipschitz continuity of operators which lead to the existence of a minimizer for these functionals. This will in turn solve IBVP-1 and IBVP-2.
\subsection{Ill-posedness of the Problems IBVP-1 and IBVP-2 }
In order to attain the compactness of the input-output operator $ \Psi$ associated with IBVP-2, we are employing more regularity on the solution to the direct problem (\ref{1-1}) as stated in Theorem \ref{8}. The proof of Theorem \ref{8} requires only $ H^2(0,T)$ regularity  of admissible source inputs, instead of $ H^3(0,T)$ regularity introduced in (\ref{1-6}), as follows
\begin{eqnarray*}
	\mathcal{G}_2=\{g\in H^2(0,T): g(0)=g^\prime(0)=0,\,\Vert g\Vert_{H^2(0,T)}\leq \bar C_g\}.
\end{eqnarray*}
But one may notice that the Lipschitz continuity of $ \Psi$ can be proved on the admissible source $ \mathcal {G}_1$ itself. 
\begin{Pro}\label{p1}
	Suppose the conditions of Theorem \ref{8} hold. Then the Neumann-Dirchlet operator $ 	\Phi :\mathcal{G}_1 \subset H^1(0,T)\mapsto L^2(0,T)$ and Neumann-Neumann operator $ \Psi:\mathcal{G}_2\subset H^2(0,T)\mapsto L^2(0,T)$ defined by (\ref{1-4}), (\ref{1-6}) respectively are compact operators. Furthermore $ \Phi$ and  $\Psi $ are Lipschitz continuous:
	\begin{eqnarray}
		\Vert \Phi (g_1)-\Phi(g_2)\Vert_{L^2(0,T)}\leq L_0 \Vert g_1 - g_2\Vert_{H^1(0,T)}, \ \ \forall \ g_1, g_2 \in \mathcal{G}_1,\label{4.10-2} \\  \label{12345}
		\Vert \Psi (g_1)-\Psi(g_2)\Vert_{L^2(0,T)}\leq L_1 \Vert g_1 - g_2\Vert_{H^1(0,T)}, \ \ \forall \ g_1, g_2 \in \mathcal{G}_1,\label{4.10-1}
	\end{eqnarray}
	with the Lipschitz constants \\ $ L_0^2= \frac{4\ell^6 (1+T) C_0^2}{9r_0^2} ,$  \
	$ L_1^2=C_{7}^2\left[1+\left(\frac {C^{\ast}\left(C_0^2+1\right)}{3\kappa_0r_0}+\frac{C_1^2}{2\rho_0}\right)(1+T)\ell^3\right],$ where \\
	$ C_{7}^2=2\ell^2 \max \left( 1,\frac{2\ell}{3}(\rho_1^2+\mu_1^2)\right)$ and $C^\ast , C_0>0$ are the constants defined in Theorem \ref{342}.
\end{Pro}
{\bf Proof.}
Using the similar arguments given in Lemma 3 and Lemma 5 of \cite{AH:OB:CS-2019}, we can show that the input-output operator $\Phi$ is compact and Lipschitz continuous with Lipschitz constant $ L_0$.
We only show that the input-output operator $ \Psi$ is compact and Lipschitz continuous with Lipschitz constant $ L_1.$

Denote by $ \{u^{(m)}(x,t)\}$, where $u^{(m)}(x,t):=u(x,t;g_m)$, the sequence of regular weak solutions  of (\ref{1-1}) corresponding to the sequence of inputs $ \{g_m\} \subset \mathcal{G}_2$, $m=1,2,...$, bounded in the norm of $ H^2(0,T).$ Then $-(r(0)u_{xx}^{(m)}(0,t)+\kappa(0)u_{xxt}^{(m)}(0,t))$ denotes the sequence of corresponding outputs. We need to prove that this sequence is relatively compact in $ L^2(0,T),$ that is, by Rellich-Kondrachov compactness theorem (see, \cite{lcevanspartial2010}, Theorem 1, Section 5.7), it is enough to show that the above sequence of outputs is bounded in the norm of $ H^1(0,T).$ To this end, we estimate the norms
$ \Vert r(0)u_{xx}^{(m)}(0,\cdot)+\kappa(0)u_{xxt}^{(m)}(0,\cdot) \Vert _{L^2(0,T)}$ and $ \Vert r(0)u_{xxt}^{(m)}(0,\cdot)+\kappa(0)u_{xxtt}^{(m)}(0,\cdot) \Vert _{L^2(0,T)}.$

Using the identity
\begin{eqnarray}\label{1234}
	r(0)u^{(m)}_{xx}(0,t) +\kappa(0)u^{(m)}_{xxt}(0,t) =\ell g_m(t)+\int_0^\ell x\left [\rho(x)u^{(m)}_{tt}+\mu(x)u^{(m)}_{t} \right ]dx,
\end{eqnarray}
which is obtained as a result of integrating (\ref{1-1}) first over the interval $(x,\ell)$ and then over the interval $(0,\ell)$, and applying   H\"older's inequality, we deduce that
\begin{eqnarray}\label{4.1.3}
	\lefteqn{	\Vert r(0)u^{(m)}_{xx}(0,.)+\kappa(0)u^{(m)}_{xxt}(0,.)\Vert^{2}_{L^2(0,T)}} \nonumber
	\\
	&&	\leq C_{7}^2 \left(\Vert g_m\Vert^2_{L^2(0,T)}+\Vert  u^{(m)}_{tt}\Vert^2_{L^2(0,T;L^2(0,\ell))}+\Vert  u^{(m)}_{t}\Vert^2_{L^2(0,T;L^2(0,\ell))}\right),
\end{eqnarray}
where $ C_{7}^2=2\ell^2 \max\{ 1,\frac{2\ell}{3}(\rho_1^2+\mu_1^2)\}.$ We infer from the estimates (\ref{1.12}) and (\ref{09}) that left-hand-side norm in (\ref{4.1.3}) is bounded in the norm of $L^2(0,T).$

Next, using the analogue
\begin{eqnarray*}
	r(0)u^{(m)}_{xxt}(0,t) +\kappa(0)u^{(m)}_{xxtt}(0,t)=\ell g_m'(t)+\int_0^\ell x\left [\rho(x)u^{(m)}_{ttt}+\mu(x)u^{(m)}_{tt} \right ]dx
\end{eqnarray*}
of identity (\ref{1234}), we deduce the estimate
\begin{eqnarray}\label{4.9}
	\lefteqn{\Vert r(0)u^{(m)}_{xxt}(0,.)+\kappa(0)u^{(m)}_{xxtt}(0,.)\Vert^{2}_{L^2(0,T)}}\nonumber \\&&
	\leq C_{7}^2 \bigg(\Vert g_m^{\prime}\Vert^2_{L^2(0,T)}+\Vert u_{ttt}^{(m)}\Vert^2_{L^2(0,T;L^2(0,\ell))}  +\Vert u_{tt}^{(m)}\Vert^2_{L^2(0,T;L^2(0,\ell))}\bigg),
\end{eqnarray}
with the same constant $C_{7}>0$. By the estimates (\ref{09}) and (\ref{4.7}), the left-hand-side norm in (\ref{4.9}) is also bounded in the norm of $L^2(0,T). $

Thus, estimates (\ref{4.1.3}) and (\ref{4.9}) imply that $ r(0)u_{xx}^{(m)}(0,t)+\kappa(0)u_{xxt}^{(m)}(0,t)$ is bounded in the norm of $ H^{1}(0,T)$. This implies that $\Psi$ is a compact operator.

In order to prove the Lipschitz continuity of $ \Psi,$  we proceed as follows. Let $ g_1, g_2 \in \mathcal{G}_1$ be the given inputs and $ u(x,t;g_1),$ $u(x,t;g_2)$ be the corresponding solutions of direct problem (\ref{1-1}).  Then $ \delta u(x,t)=u(x,t;g_1)-u(x,t;g_2)$ solves the  problem
\begin{eqnarray}   \label{1-10-0}
	\left\{ \begin{array}{lc}
		\rho(x)\delta u_{tt}+\mu(x)\delta u_t+ (r(x)\delta u_{xx})_{xx}+ (\kappa(x)\delta u_{xxt})_{xx}=0, \   (x,t)\in \Omega_T, \\ [2.5pt]
		\delta u(x,0)=0,~ \quad \; \; \ \delta u_t(x,0)=0, ~ \qquad \qquad \; \;  \ x\in (0,\ell),\\
		\delta u(0,t) \ =0,~  \quad \; \; \; \delta u_{x}(0,t)=0, ~ \qquad \qquad \quad  t\in [0,T],\\
		\left[r(x)\delta u_{xx}+\kappa(x)\delta u_{xxt}\right]_{x=\ell} \; =0, ~ \quad \qquad \ \quad \;\;\\
		\qquad -\left[\big(r(x)\delta u_{xx}+\kappa(x)\delta u_{xxt}\big)_{x}\right]_{x=\ell}=\delta g(t), \ t\in [0,T],
	\end{array}\right.
\end{eqnarray}
where $ \delta g(t)=g_1(t)-g_2(t).$
By the definition of input-output operator $ \Psi$,
\begin{eqnarray*}
	\Vert \Psi(g_1)-\Psi(g_2)\Vert^2_{L^2(0,T)}=\Vert r(0)\delta u_{xx}(0,.)+\kappa(0)\delta u_{xxt}(0,.)\Vert^2_{L^2(0,T)}.
\end{eqnarray*}
Using estimate (\ref{4.1.3})  to the solution $ \delta u(x,t)$ of problem (\ref{1-10-0}) we
arrive at the desired  estimate (\ref{12345}).  Hence the proof. \hfill$\Box$

The compactness of input-output operators $ \Phi$ and $ \Psi$ means to the ill-posedness of both inverse problems IBVP1 and IBVP2 (see,\cite{hadamard1964} and also, \cite{hasanouglu2017introduction}, Lemma 1.3.1).
\subsection{Existence and Uniqueness of Solutions to the Minimization Problems}
Using the Lipschitz continuity of the input-output operators $ \Phi$ and $ \Psi,$ we show that the existence and uniqueness of minimizer for the functionals $ \mathcal J_{1\alpha}$ and $\mathcal J_{2\alpha}$ corresponding to IBVP-1 and IBVP-2 respectively.
\begin{theorem}\label{45}
	Suppose the conditions (\ref{3-5}) hold true. Then, both minimization problems (\ref{26}) and (\ref{5.9}) have a solution on the admissible source of inputs $\mathcal{G}_1$.
\end{theorem}
{\bf Proof.}
We only prove the existence of minimizer for the functional $\mathcal J_2.$ By the same arguments, it can be proved for $\mathcal J_1$ as well.

Let $ u(x,t;g_1),$ $u(x,t;g_2)$ be the solutions of (\ref{1-1}) corresponding to the inputs $ g_1,g_2 \in \mathcal{G}_1$ respectively. Then the function $ \delta u(x,t)$ solves the problem (\ref{1-10-0}) with input $ \delta g(t)=g_1(t)-g_2(t).$
Since
\begin{eqnarray}\label{6.0}
	\Big\vert \mathcal{J}_2(g_{1})-\mathcal{J}_2(g_{2})\Big\vert^2=
	\Big\vert\sqrt{\mathcal{J}_2(g_{1})}+\sqrt{\mathcal J_2(g_{2})}\Big\vert^2\Big\vert\sqrt{\mathcal{J}_{2}(g_{1})}-\sqrt{\mathcal{J}_2(g_{2})}\Big\vert^2,
\end{eqnarray}  and appealing to the estimate (\ref{4.10-1}), we get
\begin{eqnarray}\label{7.0}
	\lefteqn{	\Big\vert\sqrt{\mathcal{J}_2(g_{1})}-\sqrt{\mathcal{J}_2(g_{2})}\Big\vert^{2}=\frac{1}{2}\Big\vert\|\Psi(g_{1})-\omega\|_{L^{2}(0,T)}-\|\Psi(g_{2})-\omega\|_{L^{2}(0,T)}\Big\vert^{2}}\nonumber\\
	&&\leq \frac{1}{2} \| \Psi(g_{1})-\Psi(g_{2})\|^{2}_{L^{2}(0,T)}\leq \frac{ L_1^2}{2} \ \Vert g_1 - g_2\Vert^2_{H^1(0,T)}, ~~~~
\end{eqnarray}
where $ L_1^2$ is defined in Proposition \ref{p1}.
Using the same step done for (\ref{4.1.3}), we can show from (\ref{1.12}) and (\ref{09}) that
\begin{eqnarray}\label{ok}
	\Vert \Psi(g_m)\Vert^2_{L^2(0,T)}  &\leq& C_{7}^2 \Big(\Vert g_m \Vert^2_{L^2(0,T)}+\Vert u_{tt}(.,.;g_m)\Vert^2_{L^2(0,T;L^2(0,\ell))} \nonumber \\ 
	&&+\Vert u_{t}(.,.;g_m)\Vert^2_{L^2(0,T;L^2(0,\ell))}\Big) 
	\leq  L_1^2 \ \Vert g_m \Vert^2_{H^1(0,T)}, 
\end{eqnarray}
$\ m=1,2.$	Applying triangle inequality and (\ref{ok}) , we obtain
\begin{eqnarray}\label{8.0}
	\Big\vert\sqrt{\mathcal{J}_2(g_{1})}+\sqrt{\mathcal{J}_2(g_{2})}\Big\vert^{2}
	&\leq& 2 \left(\| \Psi (g_{1})\|^{2}_{L^{2}(0,T)}+\| \Psi (g_{2})\|^{2}_{L^{2}(0,T)} 
	+2 \|\omega\|^{2}_{L^{2}(0,T)}\right ) \nonumber \\
	%	&\leq& 2 \left(L_{1}^2\left(\Vert g_1\Vert^2_{H^1(0,T)}+\Vert g_2\Vert^2_{H^1(0,T)}\right)+2\Vert M \Vert^2_{L^2(0,T)}\right) \nonumber \\ 
	&\leq& 4 \left( \ L_1^2 \ {C_g}+\Vert \omega \Vert^2_{L^2(0,T)}\right),
\end{eqnarray}
since $\Vert g_{m}\Vert^2_{H^1(0,T)}\leq {C_g}, \ \mbox{for} \ m=1,2.$
Consequently, (\ref{6.0}), (\ref{7.0}) and (\ref{8.0}) will lead to the estimate \begin{eqnarray*}
	\Big\vert \mathcal{J}_2(g_{1})-\mathcal{J}_2(g_{2})\Big\vert^{2}\leq 2 L_1^2 \left(  L_1^2 \ {C_g} +\Vert \omega \Vert^2_{L^2(0,T)}\right)\Vert g_1-g_2 \Vert^2_{H^1(0,T)},
\end{eqnarray*}
whence the functional $ \mathcal J_2$ is weakly lower-semi continuous on a nonempty closed convex set $\mathcal{G}_1$ (\cite{E.Zeidler}, Section 2.5, Lemma 5). Hence by the generalized Weierstrass theorem (\cite{E.Zeidler}, Theorem 2.D) the functional $\mathcal J_2(g)$ has a minimizer $ g\in \mathcal{G}_1.$  \hfill$\Box$
\begin{remark} \label{Re1}By a careful inspection, we note that the existence of solutions to IBVP-2 in Theorem {\ref{45}} is proved with inputs only in $\mathcal G_1,$ but not with  $\mathcal G_3$ defined in (\ref{1-6}). Further, the following
	Corollary \ref{Cor1} can also be proved for the Tikhonov functionals $\mathcal J_{1\alpha}$ and $ \mathcal J_{2\alpha}$ with regularizer $ \Vert g \Vert^2_{L^2(0,T)}$  instead of $\Vert g ^{\prime}\Vert^2_{L^2(0,T)}$ and $ \Vert g^{\prime \prime \prime}\Vert^2_{L^2(0,T)},$ and also can be justified that both IBVP-1, IBVP-2 have unique  solutions on $\mathcal{G}_1$.
\end{remark}
\begin{corollary} \label{Cor1}
	Assume that the conditions (\ref{3-5}) hold true. Then the regularized Tikhonov functionals $ \mathcal J_{1\alpha}(g), \ \mathcal J_{2\alpha}(g)$ defined by (\ref{54}) and (\ref{55})  have a unique minimizer on $ \mathcal{G}_1$ and $ \mathcal{G}_3,$ respectively. 
\end{corollary}
{\bf Proof.}
By Theorem {\ref{45}}, the functional $\mathcal J_{1}(g)$ is lower semi-continuous. Further, the regularized Tikhonov functional $\mathcal J_{1\alpha}(g_n)$  corresponding to $ g_{n}$ defined in (\ref{54}) satisfy $\mathcal J_{1\alpha}(g) \leq \liminf_{n\to \infty} \mathcal J_{1\alpha,n}(g),$  $\mbox{as} \ g_{n}\rightharpoonup g \ \mbox{in} \ \mathcal{G}_{1},$ 
whence $\mathcal J_{1\alpha}(g)$ is lower semi-continuous.
By the linearity of direct problem (\ref{1-1}), we have
$ u(x,t;\nu g_1+(1-\nu)g_2)=\nu u(x,t;g_1)+(1-\nu) u(x,t;g_2), \ \nu \in(0,1),$
and hence, one can get that 
%$$ \mathcal{J}_{1}(\nu g_1+(1-\nu)g_2)\leq \nu \mathcal{J}_{1}(g_1)+(1-\nu)\mathcal{J}_{1}(g_2),\ \forall \ g_1,g_2 \in \mathcal{G}_1.$$
\begin{eqnarray*}
	\mathcal{J}_{1\alpha}(\nu g_1+(1-\nu)g_2)&=& \mathcal{J}_{1}(\nu g_1+(1-\nu)g_2)+\frac{\alpha}{2}\Vert \nu g_1^{\prime}+(1-\nu)g_2^{\prime}\Vert^2_{L^2(0,T)}\nonumber \\ &<& \nu \mathcal{J}_{1\alpha}(g_1)+(1-\nu) \mathcal{J}_{1\alpha}(g_2), \forall \ g_1,g_2 \in \mathcal{G}_1, \nu \in (0,1).
\end{eqnarray*}
It shows that the functional $ \mathcal{J}_{1\alpha}(g)$ is strictly convex on $\mathcal{G}_1.$ By combining the above arguments and using the generalized Weierstrass theorem, we conclude that the functional $\mathcal J_{1\alpha}(g)$ has a unique minimizer.
By the similar arguments, we can prove that the regularized functional $\mathcal J_{2\alpha}$ defined by (\ref{55}) has a unique minimizer in admissible source inputs $ \mathcal{G}_3$. \hfill{$\Box$}
\begin{remark}
	We can directly prove the uniqueness result for the IBVP-1  as follows.
Let $ g_1,\ g_2 \in \mathcal{G}_1$ be two arbitrary given functions and $ u_{k}(x,t):= u(x,t;g_k), \ k=1,2,$ be the corresponding solutions of direct problem (\ref{1-1}).  Suppose that there exist two functions $ g_1,\ g_2\in \mathcal{G}_1,$ which are not identically zero such that $ g_1(t) \neq g_2(t),$ but the measurements $ \nu_1(t)=\nu_2(t),\ \forall \ t \in [0,T],$ where we recall that $ \nu_1(t)= u_1(\ell,t),   \nu_2(t)=u_{2}(\ell,t).$ Then the function $ z(x,t)=u_1(x,t)-u_{2}(x,t)$ solves the initial boundary value problem (\ref{1-1}) with $ g(t)$ replaced by $ \tilde{g}(t)=g_1(t)-g_{2}(t).$
Now consider the equation 
\begin{eqnarray}\label{lkj}
	\rho(x) z_{tt}+\mu(x)z_t+(r(x)z_{xx})_{xx}+(\kappa(x)z_{xxt})_{xx}=0.
\end{eqnarray}
Multiply (\ref{lkj}) by $ 2 z_{t}(x,t),$ use the formal identities (\ref{i1}), then integrating by parts, using the initial and boundary conditions, we obtain the following integral identity 
\begin{eqnarray}\label{fgd}
		\int_0^\ell \Big (\rho(x) z_{t }(t)^2+r(x)z_{xx}^2(t)\Big)dx &+&2\int_0^t \int_0^\ell \mu (x) (z_t)^2 dx d\tau  +2\int_0^t \int_0^\ell \kappa(x)(z_{xxt})^2 dx  d \tau   \nonumber \\
	&&= 2\tilde g(t)z(\ell,t)- 2 \int_0^t \tilde g^\prime(\tau)z(\ell,\tau).
\end{eqnarray}  
By the above assumption $ z(\ell,t)= u_1(\ell,t)-u_2(\ell,t)=0, t \in [0,T],$ which implies that right-hand side in (\ref{fgd}) is zero.  Since the initial and boundary conditions are homogeneous, $ z(x,t)\equiv 0, \ (x,t)\in\Omega_T.$ This contradiction completes the uniqueness of IBVP-1.

Since for the IBVP-2,  we consider measured data as bending moment $ \omega(t)$ (see, \eqref{1-5}) the identity (\ref{fgd}) cannot be used for the uniqueness of the solution.  Hence, it seems that one cannot directly get the uniqueness of the solution to the IBVP-2.  
\end{remark}
\section{Fr\'echet Differentiability of the Tikhonov Functionals}\label{s1}

In this section, we show that the  functionals $ \mathcal{J}_1(g)$ and $ \mathcal{J}_2(g)$ are Fr\'echet differentiable on the admissible sources $ \mathcal{G}_1$ and $ \mathcal{G}_3,$ respectively. The Fr\'echet derivatives are expressed in terms of the weak solutions of associated adjoint problems  with boundary data given in terms of measured data (\ref{1-2}) and (\ref{3.0}). %In order to obtain the Lipschitz continuity of Fr\'echet derivative of the functional $ \mathcal{J}_2,$  we need to assume more regularity on admissible source $\mathcal{G}_2$ that is,
%\begin{eqnarray}\label{G-3}
%\mathcal{G}_3= \{g\in \mathcal{G}_2 \cap H^0,T):g^{\prime \prime}(0)=0\},
%\end{eqnarray}
%since we are using more regularity on the solution to the direct problem (\ref{1-1}) which is given in Theorem \ref{t1}.
%But we can express the Lipschitz continuity of Frech\'et derivative of the functional $ J_1$ in terms of the solution of the adjoint problem and measured data (\ref{1-2}) without further regularity on the admissible source $ \mathcal{G}_1$.

For any $ g,\delta g \in \mathcal{G}_1$ and  $ g,\delta g \in \mathcal{G}_3$  the increments of the respective functionals $ \mathcal{J}_1(g)$ and $ \mathcal{J}_2(g)$  denoted by $ \delta \mathcal{J}_{m}(g)=\mathcal{J}_{m}(g+\delta g)-\mathcal{J}_{m}(g), \ m=1,2$ satisfies the following identities
\begin{eqnarray*}\label{56789}
	\lefteqn{\delta \mathcal{J}_1(g) = \int_0^T \left[u(\ell,t;g)-\nu(t)\right]\delta u(\ell,t) dt+\frac{1}{2}\int_0^T \delta u(\ell,t)^2 dt,} \nonumber \\
	\lefteqn{\delta \mathcal{J}_2(g)=\int_0^T  \Big(r(0)u_{xx}(0,t;g)+\kappa(0)u_{xxt}(0,t;g)+\omega(t)\Big)\Big(r(0)\delta u_{xx}(0,t)} \nonumber \\&&+\kappa(0)\delta u_{xxt}(0,t)\Big)dt +\frac{1}{2}\int_0^T \Big(r(0)\delta u_{xx}(0,t)+\kappa(0)\delta u_{xxt}(0,t)\Big)^2dt, ~~
\end{eqnarray*}
where $ \delta u(x,t):=u(x,t,g+\delta g)-u(x,t,g)$ solves the problem (\ref{1-10-0}).

The following lemma shows the representation of first integral of $ \delta \mathcal J_1(g)$ and $ \delta \mathcal J_{2}(g)$ in terms of solution of associated adjoint problems.
\begin{lemma}\label{mnbvx}
	Let assumptions (\ref{3-5}) hold true. Then the following integral relationships between the direct  and adjoint problems hold: \medskip\\
	\noindent  \textbf{(i1)} For any $ g,\ \delta g \in \mathcal{G}_1$ of IBVP-1, we have
	\begin{eqnarray}\label{90}
		\int_0^T \delta u(\ell,t) \xi(t)dt=\int_0^T \phi(\ell,t)\delta g(t)dt,
	\end{eqnarray}
	where $ \phi(x,t)$ is the solution of adjoint problem \begin{eqnarray}\label{4}
		\left\{ \begin{array}{lc}
			\rho(x)\phi_{tt}-\mu(x)\phi_t+ (r(x)\phi_{xx})_{xx}- (\kappa(x)\phi_{xxt})_{xx}=0, ~(x,t)\in \Omega_T, \\ [2.5pt]
			\phi(x,T)=0,~\phi_t(x,T)=0, ~ \qquad \qquad \quad \quad x\in (0,\ell),\\
			\phi(0,t) \ =0,~   \phi_{x}(0,t)=0, ~ \qquad \quad \quad \; \ \qquad t\in [0,T],\\
			\left[r(x)\phi_{xx}-\kappa(x)\phi_{xxt}\right]_{x=\ell} \quad =0, \\
			\qquad \left[\big(-r(x)\phi_{xx}+\kappa(x)\phi_{xxt}\big)_{x}\right]_{x=\ell}=\xi(t),  ~\quad t\in [0,T],
		\end{array}\right.
	\end{eqnarray}
	with Neumann input $ \xi\in L^2(0,T)$ and $\delta u$ is the solution of (\ref{1-10-0}). \medskip\\
	\noindent \textbf{(i2)} For any $ g, \ \delta g  \in \mathcal{G}_3$ of IBVP-2, we get
	\begin{eqnarray}\label{7}
		\int_0^T \Big(r(0)\delta u_{xx}(0,t)+\kappa(0)\delta u_{xxt}(0,t)\Big)\theta(t) dt =\int_0^T \varphi(\ell,t) \delta g(t) dt,
	\end{eqnarray}
	where $ \varphi(x,t)$ is the solution of the adjoint problem
	\begin{eqnarray}\label{2}
		\left\{ \begin{array}{lc}
			\rho(x)\varphi_{tt}-\mu(x)\varphi_t+ (r(x)\varphi_{xx})_{xx}- (\kappa(x)\varphi_{xxt})_{xx}=0, ~(x,t)\in \Omega_T, \\ [2.5pt]
			\varphi(x,T)=0,~\varphi_t(x,T)=0, ~ \qquad \qquad \quad \; x\in (0,\ell),\\
			\varphi(0,t) \ =0,~   \varphi_{x}(0,t)=\theta(t), ~ \qquad \quad \quad \; \; t\in [0,T],\\
			\left[r(x)\varphi_{xx}-\kappa(x)\varphi_{xxt}\right]_{x=\ell} \quad =0, \\
			\quad \quad \left[\big(-r(x)\varphi_{xx}+\kappa(x)\varphi_{xxt}\big)_{x}\right]_{x=\ell}=0, \  ~\quad t\in [0,T],
		\end{array}\right.
	\end{eqnarray}
	with Dirichlet input $\theta \in L^2(0,T).$
\end{lemma}
{\bf Proof.}
Multiplying the adjoint equation (\ref{4}) by $ \delta u(x, t),$ integrating over $(0,T)\times (0,\ell),$ integrating by parts and applying data values of  (\ref{1-10-0}) and (\ref{4}),  one may get
\begin{eqnarray}\label{5}
	\int_0^T \int_0^\ell \Big(\rho(x)\delta u_{t t}+\mu(x)\delta u_{t}+(r(x)\delta u_{xx})_{xx}+(\kappa(x)\delta u_{xxt})_{xx}\Big)\phi(x,t) \ dx dt\nonumber \\ [2pt]
	-\int_0^T \xi(t)\delta u(\ell,t) dt+\int_0^T \delta g(t) \phi(\ell,t) dt=0.	
	%&+&\int_0^T \Big(\Big(-r(x)\phi_{xx}(x,t)+\kappa(x)\phi_{xxt}(x,t)\Big)\delta u_x(x,t)\Big)_{x=0}^{\ell} dt\nonumber \\ [4pt]
	%&+& \int_0^T \Big(\Big(r(x)\delta u_{xx}(t)+\kappa(x)\delta u_{xxt}(t)\Big)\phi_{x}(x,t)\Big)_{x=0}^{x=\ell} dt \nonumber \\ [4pt]
\end{eqnarray}
On the other hand,   multiplying (\ref{1-10-0}) by $ \phi$ and integrating over $(0,T)\times (0,\ell),$ we notice that the first integral of (\ref{5}) is zero, which gives (\ref{90}).

Next, multiplying the first equation of (\ref{2}) by $ \delta u(x,t),$ integrate over $ (0,T)\times (0,\ell),$ and utilizing the initial and boundary conditions of (\ref{2}) and (\ref{1-10-0}), we get
\begin{eqnarray}\label{5-1-2}
	\lefteqn{\int_0^T \int_0^\ell \Big(\rho(x)\delta u_{t t}+\mu(x)\delta u_{t}+(r(x)\delta u_{xx})_{xx}+(\kappa(x)\delta u_{xxt})_{xx}\Big)\varphi(x,t) \ dx dt}\nonumber \\
	&&- \int_0^T \theta(t)\Big(r(0)\delta u_{xx}(0,t)+\kappa(0)\delta u_{xxt}(0,t)\Big) dt +\int_0^T \delta g(t)\varphi(\ell,t) dt=0. ~~~~~~~
\end{eqnarray}	
Again from equation (\ref{1-10-0}), we can conclude that the first integral of (\ref{5-1-2}) becomes zero, which leads to (\ref{7}). Hence the proof. \hfill{$\Box$}

In virtue of Theorem \ref{342}, if the arbitrary Neumann input $\xi(t)$ of adjoint problem (\ref{4}) satisfy the regularity and consistency condition $ \xi\in H^1(0,T)$, $\xi(T)=0,$ then this adjoint problem admits a unique weak solution $ \phi \in H^1(0,T;\mathcal{V}_1^2(0,\ell)),$ $ \phi_{tt}\in L^2(0,T;L^2(0,\ell))$  as the change of variable t with $ \tau=T-t$ shows. Also, $ \phi(x,t)$ satisfies the estimates (\ref{00})-(\ref{09}). Indeed, we have
\begin{eqnarray}
	\Vert \phi_{xx} \Vert^2_{L^2(0,T;L^2(0,\ell))}&\leq& \frac{4 C_{0}^2}{3r_0^2} \ (1+T) \ \ell^3\Vert \xi^{\prime}\Vert^2_{L^2(0,T)} ,\label{n41}\\ [2pt]
	\Vert \phi_{xxt}\Vert^2_{L^2(0,T;L^2(0,\ell))}&\leq& \frac{(C_0^2+1)}{3\kappa_0 \ r_0} (1+T)\ell^3 \ \Vert \xi^{\prime}\Vert^2_{L^2(0,T)},\label{bvc}\\
	\Vert \phi_{t}\Vert^2_{L^2(0,T;\mathcal{V}_1^2(0,\ell))}&\leq& \frac{ C^\ast (C_0^2+1)}{3\kappa_0 r_0} \ \ (1+T)  \ell^3 \ \Vert \xi^{\prime}\Vert^2_{L^2(0,T)},~~~~	\label{n43} \\
	\Vert \phi_{tt}\Vert^2_{L^2(0,T;L^2(0,\ell)} & \leq& \frac{C_1^2}{2\rho_0}  (1+T)\ell^3\Vert \xi^{\prime}\Vert^2_{L^2(0,T)},
\end{eqnarray}
where the constants $ C_0^2, C_1^2 , C^{\ast}, r_0,\kappa_0 $ and $\rho_0 $ are defined in Theorem \ref{342}.

The following theorem shows the necessary estimates for the weak solution $ \varphi(x,t)$ to the adjoint problem (\ref{2}).
\begin{theorem}\label{81}
	Suppose that conditions (\ref{3-5}) hold true and the Dirichlet input $ \theta(t)$ satisfies the regularity  condition $ \theta \in H^2(0,T).$ Then there exists a weak solution $ \varphi(x,t)$ of (\ref{2}) satisfying the estimates
	\begin{eqnarray}
		\Vert \varphi_t \Vert^2_{L^2(0,T;L^2(0,\ell))}&\leq& 2\left( C_{6}^2(\exp(T/\rho_0)-1)+\frac{\ell^3}{3}\right) \ G(\theta) \label{091},\\
		\Vert \varphi_{xx}\Vert^2_{L^2(0,T;L^2(0,\ell)}&\leq& \frac{T C_{6}^2}{r_0}\exp(T/\rho_{0}) \ G(\theta),\label{097} \\
		\Vert \varphi_{xxt}\Vert^2_{L^2(0,T;L^2(0,\ell))} &\leq&  \frac{C_{6}^2}{2 \kappa_0}\exp(T/\rho_{0}) \ G(\theta),\label{092}
	\end{eqnarray}
	where the constants $$ G(\theta):=\Vert \theta ^{\prime}\Vert^2_{L^2(0,T)}+\Vert \theta^{\prime \prime}\Vert^2_{L^2(0,T)},\, C_6^2= \frac{2\ell^3}{3}\left( \max (\mu_1^2,\rho_1^2)+\rho_1\max ( 1/T,\,T/3)\right).$$
\end{theorem}
{\bf Proof.}
In order to transform the adjoint problem (\ref{2}) into a problem with homogeneous boundary condition, we use the transformation (see, Appendix C.3, \cite{HRB} and also Section 3, \cite{sak2018})
$ \psi(x,t)=\varphi(x,t)-x\theta(t), \ \ (x,t)\in (0,\ell)\times[0,T).$
Then the function $ \psi(x,t)$ solves the following problem
\begin{eqnarray}\label{2.2}
	\left\{ \begin{array}{lc}
		\rho(x)\psi_{tt}-\mu(x)\psi_t+ (r(x)\psi_{xx})_{xx}- (\kappa(x)\psi_{xxt})_{xx}  \\ [2.5pt] \qquad \qquad= x\mu(x)\theta^{\prime}(t)-x\rho(x)\theta^{\prime \prime}(t),   \    (x,t)\in \Omega_T, \\ [2.5pt]
		\psi(x,T)=-x\theta(T),~\psi_t(x,T)=-x\theta^{\prime}(T), ~  \; x\in (0,\ell),\\ [2.5pt]
		\psi(0,t) \ =0,~   \psi_{x}(0,t)=0, ~ \qquad \qquad \qquad \ \; \; t\in [0,T],\\ [2.5pt]
		\left[r(x)\psi_{xx}-\kappa(x)\psi_{xxt}\right]_{x=\ell} \quad =0, \\ [2.5pt]
		\quad\quad\left[\big(-r(x)\psi_{xx}+\kappa(x)\psi_{xxt}\big)_{x}\right]_{x=\ell}=0,  \ \ t\in [0,T].
	\end{array}\right.
\end{eqnarray}
Multiply both sides of equation (\ref{2.2}) by $-2\psi_t(x,t),$ apply the identities
\begin{eqnarray*}\label{2.7}
	\left. \begin{array}{ll}
		-2\big(r(x)\psi_{xx}\big)_{xx} \psi_{t} \equiv -2 \big[(r(x)\psi_{xx})_{x} \psi_{t}-r(x)\psi_{xx} \psi_{xt}\big]_x-r(x)(\psi_{xx}^2)_t, \nonumber \\ [2.5pt]
		2\big(\kappa(x)\psi_{xxt}\big)_{xx} \psi_{t} \equiv 2 \big[(\kappa(x)\psi_{xxt})_{x} \psi_{t}+\kappa(x)\psi_{xxt} \psi_{xt}\big]_x+2\kappa(x)(\psi_{xxt}^{2}),
	\end{array} \right.
\end{eqnarray*}
and integrating by parts using the initial and boundary conditions of (\ref{2.2}), we obtain the following energy inequality
\begin{eqnarray}\label{2.6}
	\lefteqn{	\int_0^\ell \left (\rho(x) \psi_{t}(t)^2+r(x)\psi_{xx}(t)^2\right)dx +2\int_t^T \int_0^\ell \mu (x) \psi_{\tau}^2 \ dx d\tau  +2\int_t^T \int_0^\ell \kappa(x)(\psi_{xx\tau})^2 \ dx  d \tau} \nonumber \\
	&&= 2\int_t^T \int_0^\ell \Big(x\rho(x)\theta^{\prime\prime}(\tau)-x\mu(x)\theta^{\prime}(\tau)\Big)\psi_{\tau}\ dxd\tau+\int_0^\ell \rho(x)\left(x \theta^{\prime}(T)\right)^2 dx \nonumber \\
	&&\leq \int_t^T \int_0^\ell \psi_{\tau}^2 \ dx d\tau+\tilde{C}_{6}^2 \int_t^T [\theta^{\prime}(\tau)^2+\theta^{\prime \prime}(\tau)^2] d\tau + \rho_1 \frac{\ell^3}{3}\theta^{\prime}(T)^2,
\end{eqnarray}
where  $\tilde{C}_{6}^2= \frac{2\ell^3}{3} \max( \mu_1^2,\rho_1^2).$
Using the identity $ \theta^{\prime}(T)=\frac{1}{T}\int_0^T \left(t \theta^{\prime}(t)\right)_tdt,$ we get
\begin{eqnarray}\label{tre}
	\theta^{\prime}(T)^2 \leq 2 \left(\frac{1}{T} \Vert \theta ^{\prime}\Vert^2_{L^2(0,T)}+\frac{T}{3}\Vert \theta ^{\prime \prime}\Vert^2_{L^2(0,T)}\right).
\end{eqnarray}
We employ the inequality (\ref{tre}) in (\ref{2.6}) to deduce that
\begin{eqnarray}\label{2.6b}
		\int_0^\ell \Big (\rho(x) \psi_{t}(t)^2+r(x)\psi_{xx}(t)^2\Big)dx &+&2\int_t^T \int_0^\ell \mu (x) \psi_{\tau}^2 \ dx d\tau+ 2 \int_t^T \int_0^\ell \kappa(x)\psi_{xx\tau}^2 dx d\tau \nonumber\\	&&\leq \int_t^T \int_0^\ell \psi_{\tau}^2 \ dx d\tau +C_6^2\left[\Vert \theta^{\prime}\Vert^2_{L^2(0,T)}+\Vert \theta^{\prime \prime}\Vert^2_{L^2(0,T)}\right],
\end{eqnarray}
where $ C_6^2= \tilde{C}_{6}^2+\frac{2}{3}\rho_1\ell^3\max(1/T,T/3).$
Applying Gronwall's inequality, we obtain the first consequence of inequality (\ref{2.6b}) as follows
\begin{eqnarray}\label{2.03}
	\int_0^\ell \psi_{t}(t)^2 dx \leq\frac{ C_{6}^2}{\rho_0} \exp((T-t)/\rho_0)\Big[\Vert \theta ^{\prime}\Vert^2_{L^2(0,T)}+\Vert \theta^{\prime \prime}\Vert^2_{L^2(0,T)}\Big],
\end{eqnarray}
and integrate (\ref{2.03}) over $(0,T)$ to further obtain
\begin{eqnarray}\label{2.01}
	\Vert \psi_t \Vert^2_{L^2(0,T;L^2(0,\ell))}\leq C_6^2 (\exp(T/\rho_0)-1) \ G(\theta),
\end{eqnarray}
where $ G(\theta):=\Vert \theta ^{\prime}\Vert^2_{L^2(0,T)}+\Vert \theta^{\prime \prime}\Vert^2_{L^2(0,T)}.$
Substituting (\ref{2.01}) in (\ref{2.6b}), we get
\begin{eqnarray}
	\Vert \psi_{xx}\Vert^2_{L^2(0,T;L^2(0,\ell))}&\leq& \frac{T C_{6}^2}{r_0}\exp(T/\rho_{0}) \ G(\theta),\label{3.19} \\
	\Vert \psi_{xxt}\Vert^2_{L^2(0,T;L^2(0,\ell))} &\leq& \frac{C_{6}^2}{2 \kappa_0}\exp(T/\rho_{0}) \ G(\theta).\label{2.05}
\end{eqnarray}
Using the estimates (\ref{3.19}), (\ref{2.05}) and (\ref{4.8}), we get $ \psi, \ \psi_t \ \in L^2(0,T;\mathcal{V}_1^2(0,\ell)).$
Since, $\psi_t+x \theta^{\prime}(t)=\varphi_t(x,t),$  we obtain
\begin{eqnarray}
	\Vert \varphi_t \Vert^2_{L^2(0,T;L^2(0,\ell))}&\leq& 2 \Vert \psi_t \Vert^2_{L^2(0,T;L^2(0,\ell))}+\frac{2\ell^3}{3}\Vert \theta^{\prime}\Vert^2_{L^2(0,T)}, \label {2.08}
\end{eqnarray}
and $\psi_{xx}(x,t)=\varphi_{xx}(x,t),$ $ \psi_{xxt}(x,t)=\varphi_{xxt}(x,t),$
the estimates (\ref{091})-(\ref{092}) follow from (\ref{2.01})-(\ref{2.05}) and (\ref{2.08}). This completes proof. \hfill{$\Box$}

By applying the formal Lagrange multiplier method (see, \cite{FTroltzsch}, Section 3.1) for IBVP-2, we obtain the actual input $ \theta(t)$ in (\ref{2}) as follows
\begin{eqnarray}\label{098}
	\theta(t)=r(0)u_{xx}(0,t;g)+\kappa(0)u_{xxt}(0,t;g)+\omega(t), \ \ t \in [0,T].
\end{eqnarray}
Hence as a result of integral identity (\ref{7}), the following input-output relationship arises:
\begin{eqnarray*}
	\int_0^T \Big(r(0)u_{xx}(0,t)+\kappa(0)u_{xxt}(0,t)+\omega(t)\Big)\Big(r(0)\delta u_{xx}(0,t)+\kappa(0)\delta u_{xxt}(0,t)\Big)dt   =\int_0^T \varphi(\ell,t)\delta g(t) dt.
\end{eqnarray*}
Thus the variation of the functional $ \delta \mathcal J_2$ can be derived through the solution $\varphi(x,t)$ of the adjoint problem (\ref{2}) with the input (\ref{098}) as follows:
\begin{eqnarray}\label{9.6}
	\delta \mathcal{J}_2(g)\!\!=\!\!\int_0^T \varphi(\ell,t)\delta g(t) dt+\frac{1}{2}\int_0^T \left(r(0)\delta u_{xx}(0,t)+\kappa(0)\delta u_{xxt}(0,t)\right)^2 dt.
\end{eqnarray}
Similarly, we consider the  adjoint problem (\ref{4}) with Dirichlet input
\begin{eqnarray}\label{iklj}
	\xi(t)= u(\ell,t;g)-\nu(t), \ t\in[0,T].
\end{eqnarray}
As in the case of IBVP-2, we also obtain that
$
\int_0^T \left[u(\ell,t;g)-\nu(t)\right]\delta u(\ell,t)=\int_0^T \phi(\ell,t)\delta g(t)dt,
$
and
\begin{eqnarray}\label{5-1-3}
	\delta \mathcal{J}_1(g)=\int_0^T \phi(\ell,t)\delta g(t)dt+ \frac{1}{2}\int_0^T \delta u(\ell,t)^2 dt .
\end{eqnarray}

Let us now justify the substitutions (\ref{iklj}) and (\ref{098}) in the context of solutions of (\ref{4}) and (\ref{2}), respectively. By invoking the theory developed in (\cite{Baysal:Hasanov:2019}) to the backward problem (\ref{2.2}) and Theorem \ref{81}, it follows that the weak solution $ \varphi \in H^1([0,T];\mathcal{V}_1^2(0,\ell))$ of problem (\ref{2}) exists, if the input $\theta(t)$ belongs to $H^2(0, T)$. Furthermore, by Theorem {\ref{342}}, for the existence of the weak solution $\phi \in H^1([0,T];\mathcal{V}_1^2(0,\ell))$ of problem (\ref{4}), the Neumann input $\xi(t)$ should satisfy
the regularity condition $\xi \in H^1(0,T)$. In view of the substitutions (\ref{iklj}) and (\ref{098}), we infer that  the measured outputs $\nu(t)$, $\omega(t)$ and the outputs $u(\ell,t;g)$,  $-r(0)u_{xx}(0,t;g)-\kappa(0)u_{xxt}(0,t;g)$ should obey the following regularity conditions:
\begin{eqnarray}
	&& \nu \in H^1(0,T),~\omega\in H^2(0,T), \label{5-mout} \\ [4pt]
	&&u(\ell,\cdot) \in H^1(0,T),~-r(0)u_{xx}(0,\cdot)-\kappa(0)u_{xxt}(0,\cdot)\in H^2(0,T). \label{5-out}
\end{eqnarray}

The conditions (\ref{5-mout}) mean that the measured outputs $\nu(t)$ and $\omega(t)$ should be more regular, although both of them originally belong to the class $L^2(0, T)$. In particular, the Neumann measured output $\omega(t)$ requires more regularity than that of the Dirichlet measured output $\nu(t)$. Furthermore, evidently for the weak solution of the direct problem (\ref{1-1}), the first condition of (\ref{5-out}) holds.  The identity
\begin{eqnarray}\label{ident-1}
	r(0)u_{xxtt}(0,t) +\kappa(0)u_{xxttt}(0,t)  =\ell g''(t)+\int_0^\ell x\left [\rho(x)u_{tttt}+\mu(x)u_{ttt} \right ]dx, \quad \forall t\in[0,T],
\end{eqnarray}
for the Neumann output and Theorem \ref{t1}  show that  the second condition of (\ref{5-out}) is also satisfied for the regular weak solution with improved regularity.
\begin{theorem}\label{tr} Let the conditions ({\ref{3-5}}) hold true.\\
	\noindent \textbf{(i1)}	Assume that the Dirichlet measured output $\nu(t)$ satisfies the regularity condition $\nu \in H^1(0,T)$. Then the Tikhonov functional $\mathcal J_1(g)$ corresponding to the problem IBVP-1 is Fr\'echet differentiable on the set of admissible sources $ \mathcal{G}_1$. Furthermore, for the Fr\'echet derivative of this functional at $ g\in\mathcal{G}_1,$ the following gradient formula holds:
	\begin{eqnarray}\label {6.8}
		\mathcal{J}_1^{\prime}(g)(t)=\phi(\ell,t;g,\xi), \, t\in(0,T),
	\end{eqnarray}
	where $\phi(x,t;g)$ is the weak solutions of the adjoint problem (\ref{4}) with the input  $\xi(t)=u(\ell,t;g)-\nu(t)$.
	
	\noindent \textbf{(i2)} Assume that the conditions of Theorem \ref{t1} hold and $g\in \mathcal{G}_3$. Suppose in addition,  the Neumann measured output $\omega(t)$ satisfies the regularity condition $\omega \in H^2(0,T)$. Then the Tikhonov functional $\mathcal J_2(g)$ corresponding to the problem IBVP-2 is Fr\'echet differentiable on $ \mathcal{G}_3$. Moreover, the gradient formula
	\begin{eqnarray}\label {6.9}
		\mathcal{J}_2^{\prime}(g)(t)=\varphi(\ell,t;g,\theta), \, t\in(0,T)
	\end{eqnarray}
	holds through the weak solution $\varphi(x,t;g,\theta)$ of the adjoint problem (\ref{2}) with the input  $\theta(t)=r(0)u_{xx}(0,t;g)+\kappa(0)u_{xxt}(0,t;g)+\omega(t)$.
\end{theorem}
{\bf Proof.}
By employing the inequality (\ref{1-10}) and estimate (\ref{1.2}) to the solution $\delta u(x,t)$ of problem (\ref{1-10-0}), we obtain from (\ref{5-1-3}) that
\begin{eqnarray}\label{ty}
	\Big\vert \delta \mathcal{J}_1(g)-\int_0^T \phi(\ell,t) \delta g(t) dt \Big\vert&=& \frac{1}{2} \Vert \delta u (\ell,.)\Vert^2_{L^2(0,T) } \leq \frac{\ell^3}{6} \Vert \delta u_{xx}\Vert^2_{L^2(0,T;L^2(0,\ell))} 
	\nonumber \\ &&\leq \frac{2 \ell^6 (1+T)C_0^2}{9r_0^2} \ \Vert \delta g \Vert^2_{H^1(0,T)},
\end{eqnarray}
where $ \delta g(t)=g_{1}(t)-g_{2}(t).$
As a consequence, we have
$$ \frac{	\Big\vert \delta \mathcal{J}_1(g)-\int_0^T \phi(\ell,t) \delta g(t) dt \Big\vert}{\Vert \delta g \Vert_{H^1(0,T)}}\to 0  \ \ \  \mbox{as} \ \ \|\delta g\|_{H^{1}(0,T)} \to 0^+.$$
This means the Fr\'echet differentiability of the functional $\mathcal J_1$.

Next, we consider formula (\ref{9.6}) for the variation of the functional $\mathcal J_2$. Using estimate (\ref{4.1.3})  to the solution $\delta u(x,t)$ of  (\ref{2.2}), and then estimates (\ref{1.12}) and (\ref{09}), we get
\begin{eqnarray}\label{tu}
	\lefteqn{\Big\vert \delta \mathcal{J}_2(g)-\int_0^T \varphi(\ell,t) \delta g(t) dt \Big\vert=\frac{1}{2}\Vert r(0)\delta u_{xx}(0,.)+\kappa(0)\delta u_{xxt}(0,.)\Vert^2_{L^2(0,T)}}\nonumber \\ &\leq& \frac{C_{7}^2}{2} \Big(\Vert \delta u_{tt}\Vert^2_{L^2(0,T;L^2(0,\ell))}+\Vert \delta u_{t}\Vert^2_{L^2(0,T;L^2(0,\ell))} +\Vert \delta g\Vert^2_{L^2(0,T)}\Big)
	\nonumber \\ &	\leq& C_{8}^2 \Vert \delta g \Vert^2_{H^1(0,T)},
\end{eqnarray}
where $ C_{8}^2=\frac{C_{7}^2}{2}\left[1+\left(\frac{C_1^2}{2\rho_0}+\frac{C^\ast \left(C_0^2+1\right)}{3r_0\kappa_0}\right) (1+T)\ell^3\right].$  This shows that  $$ \frac{	\vline \delta \mathcal{J}_2(g)-\int_0^T \varphi(\ell,t) \delta g(t) dt \vline}{\Vert \delta g \Vert_{H^1(0,T)}}\to 0  \ \ \  \mbox{as} \ \ \|\delta g\|_{H^{1}(0,T)} \to 0^+.$$
Hence, the definition of the Fr\'echet derivative gives the formulas (\ref{6.8}) and (\ref{6.9}). \hfill{$\Box$}

\begin{remark}\label{r2}
	The gradient formulas (\ref{6.8}) and (\ref{6.9}) show that there is no  need for the
	weak solutions $\phi(x,t)$ and $\varphi(x,t)$ of adjoint problems (\ref{4}) and (\ref{2}) with inputs
	(\ref{iklj}) and (\ref{098}), respectively. Namely, these solutions need to satisfy only the
	conditions $\phi_x, \varphi_x  \in L^2(0,T;L^2(0,\ell)).$ By introducing a weaker solution, as in   (\cite{Hao}) and (\cite{hasanouglu2019}), the above conditions (\ref{5-mout}) and (\ref{5-out}) can be
	weakened.
\end{remark}

\begin{corollary}
	Suppose the conditions of Theorem \ref{tr} hold true. Then the regularized Tikhonov functionals  $ \mathcal J_{1\alpha}(g)$, $\mathcal J_{2\alpha}(g)$ defined in (\ref{54}) and (\ref{55}) are Fr\'echet differentiable on $ \mathcal{G}_1$ and $\mathcal{G}_3,$ respectively. The Fr\'echet derivatives are given by
	\begin{eqnarray}
		\mathcal{J}_{1\alpha}^{\prime}(g)(t)&=&\phi(\ell,t;g,\xi)+\alpha g^{\prime}(t), \ \forall g\in \mathcal{G}_{1}\label{3452}\\
		\mathcal{J}_{2\alpha}^{\prime}(g)(t)&=&\varphi(\ell,t;g,\theta)+\alpha g^{\prime \prime \prime}(t), \  \forall g \in \mathcal{G}_3,\label{3478}
	\end{eqnarray}
	where $ \phi, \,\varphi  \in H^1(0,T;\mathcal{V}_1^2(0,\ell))$ are the weak solutions of the adjoint problems (\ref{4}) and (\ref{2})  with boundary data $\xi$ and $\theta$ as in Theorem \ref{tr}, respectively.
\end{corollary}
{\bf Proof. }
For any $ g,\ \delta g \in \mathcal{G}_1$ and $ g,\ \delta g \in \mathcal{G}_3,$ the increment corresponding to the functionals  $\mathcal J_{m\alpha}(g), \ m=1,2$ are given by
\begin{eqnarray*}
	\delta \mathcal{J}_{1\alpha}(g)&= &\int_0^T \phi(\ell,t)\delta g(t)dt+ \alpha\int_0^T g^{\prime}(t) \delta g^{\prime}(t) dt+ \frac{\alpha}{2} \int_0 ^T (\delta g^{\prime}(t) )^2dt +\frac{1}{2}\int_0^T \delta u(\ell,t)^2 dt,\ \ \forall g, \delta g \in \mathcal{G}_1,\label{2345}\\
	\delta \mathcal{J}_{2\alpha}(g)&=& \int_0^T \varphi(\ell,t)\delta g(t) dt +\alpha \int_0^T  g^{\prime \prime \prime}(t) \delta g^{\prime \prime \prime}(t) dt+\frac{\alpha}{2}\int_0^T(\delta g^{\prime \prime \prime}(t))^2 dt \nonumber \\ &&+\frac{1}{2}\int_0^T \Big(r(0)\delta u_{xx}(0,t)+\kappa(0)\delta u_{xxt}(0,t)\Big)^2 dt, \ \ \forall g, \delta g \in \mathcal{G}_3.
\end{eqnarray*}
By doing calculations similar to (\ref{ty}) and (\ref{tu}) of Theorem \ref{tr}, we can show that the last two integrals of $ \delta \mathcal J_{m\alpha},\ m=1,2$ are of the orders $\mathcal{O} \left(\Vert\delta g\Vert^2_{H^1(0,T)}\right)$ and $\mathcal{O} \left(\Vert\delta g\Vert^2_{H^3(0,T)}\right)$, respectively. From the definition of Fr\'echet derivative, we obtain the desired results (\ref{3452}) and (\ref{3478}). 	\hfill{$\Box$}

\section{Monotonicity of the Gradient Algorithm} \label{s55}
The Lipschitz continuity of the Fr\'echet derivatives of functionals $ \mathcal{J}_1, \mathcal{J}_2 $ has an important advantage when applying gradient-based methods to solve an inverse problem. In particular, in the case of gradient type algorithms such as Landweber iteration algorithm  $ g^{(n+1)}(x)=g^{(n)}(x)-\gamma_{n}\mathcal{J}^{\prime}(g^{(n)}(x)), \ n=0,1,2,...,$ or conjugate gradient algorithm applied to solve inverse problems, we may have trouble in predicting the relaxation parameter $ \gamma_n > 0.$ Using the Lipschitz constants associated with the Lipschitz continuity of $ \mathcal{J}_1^{\prime}, \mathcal{J}_2^{\prime},$ the relaxation parameter can be calculated and that can be used to discuss the convergence of the iterative scheme as well.
The following result shows that $\mathcal{J}_1^{\prime}$, $ \mathcal{J}_2^{\prime}$ are Lipschitz contionus on $ \mathcal{G}_1$ and $\mathcal{G}_3$ respectively.
\begin{Pro}\label{3}
	Let the conditions of Theorem \ref{tr} hold true. Then the Fr\'echet gradients of the functionals $ \mathcal{J}_{i}^{\prime}(g),\ i=1,2,$ defined by (\ref{6.8}) and (\ref{6.9}) are Lipschitz continuous. Moreover,
	\begin{eqnarray}
		\Vert \mathcal{J}_1^{\prime}(g+\delta g)-\mathcal{J}_1^{\prime}(g)\Vert_{L^2(0,T)}\leq L_2 \Vert \delta g \Vert_{H^1(0,T)}, \ g,\delta g \in \mathcal{G}_1,\label{48} \\
		\Vert \mathcal{J}_2^{\prime}(g+\delta g)-\mathcal{J}_2^{\prime}(g)\Vert_{L^2(0,T)}\leq L_3 \Vert \delta g \Vert_{H^3(0,T)}, \ g,\delta g \in \mathcal{G}_3,\label{49}
	\end{eqnarray}
	where the Lipschitz constants
	\begin{eqnarray*}
		&& L_2^2=  \frac {C_0^2}{\kappa_0 r_0}(C_0^2+1) \left(\frac{2\ell^6}{9 r_0} \  (1+T) \right)^2,\nonumber \\ &&  L_3^2= \Big(\frac{C_{7}^2\ell^3 T C_6^2}{3 r_0}\exp(T/\rho_0)\Big) \left(1+\frac{1}{2\rho_0}\left[3C_5^2 \exp(C_5^2 T)+C_1^2 \ell^3(1+T) \right]\right),
	\end{eqnarray*}
	where $ C_{5}, C_{6}, C_{7}>0$ are the constants introduced in Theorem \ref{8}, Theorem \ref{81} and Proposition \ref{p1}, respectively.
\end{Pro}

{\bf Proof.}
For any $ g,\delta g \in \mathcal{G}_1,$ from the Fr\'echet derivative (\ref{6.8}) it is clear that $ \Vert\mathcal{J}_1^{\prime}(g+\delta g)-\mathcal{J}_1^{\prime}(g)\Vert_{L^2(0,T)}^{2} = \int_0^T \delta \phi(\ell,t)^2 dt,$ where $ \delta \phi $ is the solution of (\ref{4}) with data $ \delta\xi(t)= \delta u(\ell,t)=u(\ell,t;g+\delta g)-u(\ell,t;g).$ Applying the trace estimate (\ref{1-10}) which holds for $ \delta \phi(\ell,t),$ and the estimates (\ref{n41}), (\ref{2.98}), we obtain
\begin{eqnarray}
	\Vert \mathcal{J}_1^{\prime}(g+\delta g)-\mathcal{J}_1^{\prime}(g)\Vert_{L^2(0,T)}^{2}&\leq& \frac{\ell^3}{3}\Vert \delta \phi_{xx}\Vert^2_{L^2(0,T;L^2(0,\ell))} 
	\leq \frac{4\ell^6}{9 r_0^2}  C_0^2 (1+T) \Vert \delta \xi^{\prime} \Vert^2_{L^2(0,T)} \nonumber \\ 
	&\leq& \frac{4\ell^9 }{27 r_{0}^2} C_0^2 (1+T)   \Vert \delta u _{xxt}\Vert^2_{L^2(0,T;L^2(0,\ell))}\leq  L_2^2 \Vert \delta g \Vert^2_{H^1(0,T)},
\end{eqnarray}
%In order to apply the estimate (\ref{n41}) with $ \xi(t)=\delta u(\ell,t),$ we need to show that $\delta u(\ell,t) \in H^1(0,T).$ This follows from inequalities  (\ref{1-10}) and (\ref{mnbv}), with  (\ref{n4}) and (\ref{1.12}).  Hence we can apply the estimate  (\ref{n41}) and we get,
%since we also used $ \Vert \delta \xi^{\prime}\Vert^2_{L^2(0,T)}=\Vert \delta u_t(\ell,.)\Vert^2_{L^2(0,T)}\leq \frac{\ell^3}{3}  \Vert \delta u _{xxt}\Vert^2_{L^2(0,T;L^2(0,\ell))},$	
where $L_2^2= \frac {C_0^2}{\kappa_0 r_0}(C_0^2+1) \left(\frac{2\ell^6}{9 r_0} \  (1+T)  \right)^2.$ 
Next we prove that $\mathcal J_2^{\prime}(g)$ is Lipschitz continuous. The gradient formula  (\ref{6.9}) and the trace inequality (\ref{1-10})  lead to
\begin{eqnarray*}
	\Vert \mathcal{J}_2^{\prime}(g+\delta g)-\mathcal{J}_2^{\prime}(g)\Vert_{L^2(0,T)}^{2} = \Vert \delta \varphi(\ell,.)\Vert^2_{L^2(0,T)}\leq \frac{\ell^3}{3} \Vert \delta \varphi_{xx}\Vert^2_{L^2(0,T;L^2(0,\ell))},
\end{eqnarray*}
where $ \delta \varphi (x,t)$ is the solution of (\ref{2}) with the boundary data $ \delta \varphi_{x}(0,t)=\delta\theta(t) =\left( r(0) \delta u_{xx}(0,t) \\ +\kappa(0)\delta u_{xxt}(0,t)\right).$
By employing the estimate (\ref{097}), one can get
\begin{eqnarray}
	\Vert \mathcal{J}_2^{\prime}(g+\delta g)-\mathcal{J}_2^{\prime}(g)\Vert_{L^2(0,T)}^{2} &\leq&
		\frac{\ell^3  T  C_6^2}{3 r_0}\exp(T/\rho_0) \Big( \Vert r(0)\delta u_{xxt}(0,.)+\kappa(0)\delta u_{xxtt}(0,.)\Vert^2_{L^2(0,T)}\nonumber \\ 
	&&+\Vert r(0)\delta u_{xxtt}(0,.)+\kappa(0)\delta u_{xxttt}(0,.)\Vert^2_{L^2(0,T)} \Big)~~.
\end{eqnarray}
From the identity (\ref{ident-1}), it holds that
\begin{eqnarray} \label{076}
	\lefteqn{\Vert\left(r(0)\delta u_{xx}(0,.;g)+\kappa(0)\delta u_{xxt}(0,.;g)\right)_{tt}\Vert^2_{L^2(0,T)}} \nonumber \\ [4pt]&\leq &C_{7}^2 \left(\Vert \delta u_{tttt}\Vert^2_{L^2(0,T;L^2(0,\ell))}+\Vert\delta u_{ttt}\Vert^2_{L^2(0,T;L^2(0,\ell))}+\Vert \delta g^{\prime \prime}\Vert^2_{L^2(0,T)}\right),
\end{eqnarray}
where $C_7$ is the constant defined in \ref{4.1.3}, and coupling with  (\ref{4.9}), we arrive at 
\begin{eqnarray}
	\lefteqn{	\Vert \mathcal{J}_2^{\prime}(g+\delta g)-\mathcal{J}_2^{\prime}(g)\Vert_{L^2(0,T)}^{2}} \nonumber \\ && \leq \frac{ C_{7}^2 \ell^3  T  C_6^2}{3 r_0}\exp(T/\rho_0) \Big(\Vert \delta g^{\prime \prime} \Vert^2_{L^2(0,T)}+ \Vert \delta g^{\prime}\Vert^2_{L^2(0,T)} +\Vert \delta u_{tttt}\Vert ^2_{L^2(0,T;L^2(0,\ell))}\nonumber \\ 
	&&+2\Vert \delta u_{ttt}\Vert ^2_{L^2(0,T;L^2(0,\ell))}+ \Vert \delta u_{tt}\Vert ^2_{L^2(0,T;L^2(0,\ell))}\Big).
\end{eqnarray}
By the estimates (\ref{09}), (\ref{4.1}) and (\ref{4.7}), we obtain the desired result (\ref{49}). \hfill{$\Box$}

Next, we discuss the convergence of Landweber iterative  scheme.
The sequence of iterations $\{g^{(n)}\}\subset \mathcal{G}_1$ of IBVP-1 and $\{g^{(n)}\} \subset \mathcal{G}_3$ of IBVP-2 are  defined by
\begin{eqnarray}\label{87}
	g^{(n+1)}(t)=g^{(n)}(t)-\gamma_{n}\mathcal{J}_m^{\prime}(g^{(n)}(t)), \ m=1,2, \ n=0,1,2,...,
\end{eqnarray}
where iteration parameter $ \gamma_n$ is given by the minimum problem
$$ f_{n}(\gamma_n):=\inf_{\gamma\geq 0}f_{n}(\gamma), \ \ f_{n}(\gamma):= \mathcal{J}(g^{(n)}-\gamma \mathcal{J}_{m}^{\prime}(g^{(n)}))(t).$$
\begin{Pro}
	Assume that the conditions of Proposition \ref{3} hold true and let $ g^{(n)}$ be the iteration defined by (\ref{87}) with $ \gamma_n=\gamma >0. $ Then the following inequalities hold
	\begin{eqnarray*}
		\mathcal{J}_{1}(g^{(n)})-\mathcal{J}_{1}(g^{(n+1)})& \geq& \frac{1}{2L_2}\Vert \mathcal{J}_{1}^{\prime}(g^{(n)})\Vert^2_{L^2(0,T)},\nonumber \\
		\mathcal{J}_{2}(g^{(n)})-\mathcal J_{2}(g^{(n+1)}) &\geq& \frac{1}{2L_3}\Vert \mathcal{J}_{2}^{\prime}(g^{(n)})\Vert^2_{L^2(0,T)}, \ \ \forall \ n=0,1,2,... ,
	\end{eqnarray*}
	where $ L_2$ and $ \ L_3$ are the Lipschitz constants defined in Proposition \ref{3}.
	Moreover, the sequence $\mathcal{J}_{m}(g^{(n)}),\ m=1,2$ is monotone decreasing convergent sequence with  $\lim_{n\to\infty} \Vert \mathcal{J}_{m}^{\prime}(g^{(n)})\Vert_{L^2(0,T)} = 0,  \ m=1,2.$
\end{Pro}
The proof follows from the similar arguments of Lemma 4.3 and Corollary 4.1 of \cite{hasanov2007}. Consequently, let us set $ \mathcal{J}_{1}^{\ast}=\mathcal{J}_{1}(g^{\ast})=\lim_{n\to \infty}\mathcal{J}_{1}(g^{(n)})$ be the limit of the sequence $ \mathcal{J}_{1}(g^{(n)}).$ It is evident that the sequence of iterations $\{ g^{(n)}\}\subset\mathcal{G}_1$ of IBVP-1 weakly converges to $ g^{\ast}$ in $ L^2(0,T)$  . Similar, conclusions hold for IBVP-2 as well.
\section{Stability Estimates by Variational Methods}\label{s2}
This section establishes a variational inequality, which has to be satisfied by an optimal solution of the minimization problems (\ref{54}) and (\ref{55}).
This variational inequality is the key ingredient in deriving the stability estimates for the inverse problems. The stability estimates for IBVP-1 and IBVP-2 are obtained through the regular solutions established in Theorems \ref{342},\ref{8} and \ref{t1}  under a suitable smoothness of the boundary data $g(t).$ This forces us  to introduce the more regularized Tikhonov functionals $ \mathcal{J}_{1\alpha}(g)$ and $ \mathcal{J}_{2\alpha}(g)$ as in (\ref{54}) and (\ref{55}) with inputs  in $\mathcal{G}_1$ and $\mathcal{G}_3$, respectively.
\begin{Pro}\label{J}
	Let $(\bar{u}(\ell,.),\bar{g})$ and $(-(r(0)\bar{u}_{xx}(0,.)+\kappa(0)\bar{u}_{xxt}(0,.)),\bar{g})$ be the solutions of IBVP-1 and IBVP-2 respectively. Then for the problem IBVP-1,  the following variational inequality holds
	\begin{eqnarray}\label{546}
		\int_{0}^{T} \left(g_{ \alpha}(t)-\bar{g}(t)\right)\phi(\ell,t;\bar g) dt+\alpha \int_0^T \bar{g}^{\prime}(t)\left(g_{\alpha}^{\prime}(t)-\bar{g}^\prime(t)\right)  dt \geq 0,
	\end{eqnarray}
	$\forall  g_{\alpha}\in \mathcal{G}_1,$ while for the case of IBVP-2, it holds that
	\begin{eqnarray} \label{547}
		\int_{0}^{T} \left(g_{ \alpha}(t)-\bar{g}(t)\right)\varphi(\ell,t;\bar g)  dt  +\alpha \int_0^T \bar{g}^{\prime\prime \prime}(t)\left(g_{\alpha}^{\prime \prime\prime}(t)-\bar{g}^{\prime \prime \prime}(t)\right)  dt \geq 0,  
	\end{eqnarray}
	$\forall g_{\alpha}\in \mathcal{G}_3,$	where $\phi$ and $ \varphi$ are the weak solutions of the adjoint problems (\ref{4}) and (\ref{2}) with data (\ref{iklj}) and (\ref{098}), respectively.	
\end{Pro}
{\bf Proof.}
For any $ 0\leq \gamma\leq1,$ we choose an arbitrary element $ g_{ \alpha}\in\mathcal{G}_1$ such that  $ g_{\gamma}=\bar{g}+\gamma(g_{ \alpha}-\bar{g}) \in \mathcal{G}_1.$ The regularized Tikhonov functional corresponding to $(u_{\gamma}(\ell,.),g_{\gamma})$ is given by
\begin{eqnarray*}
	\mathcal{J}_{1\alpha}(g_{\gamma})=\frac{1}{2}\int_0^T \left(u_{\gamma}(\ell,t;g_{\gamma})-\nu(t)\right)^2 dt +\frac{\alpha}{2}\int_0^T \left(g_{\gamma}^{\prime}(t)\right)^2 dt.
\end{eqnarray*}
Since the functional $ J_{1\alpha}(g_{\gamma})$ is Fr\'echet differentiable at $ g_{\gamma},$ we have
\begin{eqnarray} \label{jbb} \frac{d}{d\gamma}\Big(\mathcal{J}_{1\alpha}\left(\bar{g}+\gamma(g_{\alpha}-\bar{g})\right)\Big)\Big\vert_{\gamma=0}=\int_0^T \Big(u_{\gamma}(\ell,t;g_\gamma)-\nu(t)\Big)\frac{\partial u_{\gamma}}{\partial \gamma}\vert_{\gamma=0} \ dt +\alpha \int_0^T \bar{g}^{\prime}(t)\left(g_{ \alpha}^{\prime}(t)-\bar{g}^{\prime}(t)\right)dt.
\end{eqnarray}
Considering the system (\ref{1-1}) corresponding to the data $ g_{\gamma}$ and setting	$\eta=\frac{\partial u_{\gamma}}{\partial \gamma }\vert_{\gamma=0},$ we see that $ \eta(x,t)$ satisfies the system
\begin{eqnarray}\label{5456}
	\left\{ \begin{array}{lc}
		\rho(x)\eta_{tt}+\mu(x)\eta_t+ (r(x)\eta_{xx})_{xx}+ (\kappa(x)\eta_{xxt})_{xx}=0, ~(x,t)\in \Omega_T, \\ [2.5pt]
		\eta(x,0)=0,~\eta_t(x,0)=0, ~ \qquad \; x\in (0,\ell),\\ [2.5pt]
		\eta(0,t) \ =0,~    \; \eta_{x}(0,t)=0, ~  \qquad \; t\in [0,T],\\ [2.5pt]
		\left[r(x)\eta_{xx}+\kappa(x)\eta_{xxt}\right]_{x=\ell} \; =0, \\ [2.5pt]
		\quad-\left[\big(r(x)\eta_{xx}+\kappa(x)\eta_{xxt}\big)_{x}\right]_{x=\ell}=g_{\alpha}(t)-\bar{g}(t),  \ \ t\in [0,T].
	\end{array}\right.
\end{eqnarray}
Since $ \bar{g} $ is the optimal solution, we obtain
$$ \frac{d}{d\gamma}\Big(\mathcal{J}_{1\alpha}(\bar{g}+\gamma (g_{ \alpha}-\bar{g}))\Big)\Big\vert_{\gamma=0} \geq 0, \ \ \forall \ g_{ \alpha}\in \mathcal{G}_1.$$
Therefore, from (\ref{jbb}) we have
$$ \int_{0}^{T} \left[\bar{u}(\ell,t;\bar{g})-\nu(t)\right]\eta(\ell,t) dt +\alpha \int_0^T \bar{g}^{\prime}(t)\left(g^{\prime}_{ \alpha}(t)-\bar{g}^{\prime}(t)\right)dt \geq 0, \; \; \forall \ g_{ \alpha}\in\mathcal{G}_1.$$
In virtue of the relationship $\left[\bar{u}(\ell,t;\bar{g})-\nu(t)\right] =\left( -r(x)\phi_{xx}+\kappa(x)\phi_{xxt}\right)_{x}\vert_{x=\ell},$ between the measured data and the adjoint solution $ \phi$ of (\ref{4}), one can rewrite
\begin{eqnarray}\label{543}
	\int_{0}^{T} \left( -r(x)\phi_{xx}+\kappa(x)\phi_{xxt}\right)_{x}\vert_{x=\ell} \ \eta(\ell,t) dt +\alpha \int_0^T \bar{g}^{\prime}(t)\left(g^{\prime}_{ \alpha}(t)-\bar{g}^{\prime}(t)\right)dt \geq 0.~~~~~
\end{eqnarray}
In order to express the first integral of (\ref{543}) solely interms of solution of adjoint system, we multiply equation (\ref{4}) by $ \eta(x,t),$ integrating by parts and apply initial and boundary conditions, we get
\begin{eqnarray*}
	&-&\int_0^T \left(-r(x)\phi_{xx}+\kappa(x)\phi_{xxt}\right)_{x}\vert_{x=\ell}\ \eta(\ell,t)dt \nonumber \\  &+&\int_{0}^{\ell} \int_{0}^T \Big(\rho(x)\eta_{tt}+\mu(x)\eta_{t}+(r(x)\eta_{xx})_{xx}+(\kappa(x)\eta_{xxt})_{xx}\Big)\phi(x,t) dt  dx  + \int_0^T \left(g_{ \alpha}(t)-\bar{g}(t)\right)\phi(\ell,t) dt=0.
\end{eqnarray*}
Multiplying (\ref{5456}) by $ \phi(x,t)$, integrating over $ (0,\ell)\times (0,T),$ and using it in the previous equation, we get
$$ \int_0^T \left(-r(x)\phi_{xx}+\kappa(x)\phi_{xxt}\right)_{x}\vert_{x=\ell}\ \eta(\ell,t)\ dt= \int_0^T \left(g_{ \alpha}(t)-\bar{g}(t)\right)\phi(\ell,t) dt.$$
Substitution of this identity in (\ref{543}) leads to the desired inequality (\ref{546}).

By repeating the calculation for $ \mathcal{J}_{2\alpha}(g_\gamma)$  with the adjoint problem (\ref{2}) and the input  (\ref{098}), one can obtain the variational inequality (\ref{547}) for IBVP-2. This completes the proof.
\hfill{$\Box$}

Next, we have the following stability estimate for the IBVP-1 in terms of the measured data. We obtain a lower bound for the internal damping coefficient $ \kappa(x)$ which is sufficient to obtain a Lipschitz type stability estimate for the shear force $ g(t).$
\begin{theorem}\label{mnk}
	Suppose the assumptions (\ref{3-5}) hold true. Let $ g_{\alpha}, \ \widehat{g}_{\alpha} \in \mathcal{G}_1$ are unique minimizers of the regularized Tikhonov functional $\mathcal  J_{1\alpha}$ defined by (\ref{54}) corresponding to the measured outputs $ \nu,\ \widehat{\nu} \in H^1(0,T),$ respectively.
	Suppose the internal damping coefficient $ \kappa(x)$ satisfies the condition
	\begin{eqnarray}\label{643}
		\kappa(x)\geq \Big(\frac{\sqrt 2T^2 \ell^6\exp(T)(1+T)}{9r_0}\Big)\alpha^{-1}:=\kappa_0.
	\end{eqnarray}
	Then the following stability estimate holds:
	\begin{eqnarray}\label{646}
		\Vert \widehat{g}_{\alpha}-g_{ \alpha}\Vert^2_{L^2(0,T)}\leq C_{ST}\Vert \widehat{\nu}^{\prime}-\nu^{\prime}\Vert^2_{L^2(0,T)},
	\end{eqnarray}
	where $ C_{ST}=\frac{9 r_0\kappa_0 T^2}{\ell^6  \exp(T)(1+T)}.$
\end{theorem}
{\bf Proof.}
Take $  \bar{g}(t)=\widehat{g}_{ \alpha}(t)$ in the variational inequality (\ref{546}) to get
\begin{eqnarray}\label{bn}
	\int_{0}^{T} \left(g_{ \alpha}(t)-\widehat{g}_{ \alpha}(t)\right)\phi(\ell,t;\widehat{g}_{ \alpha}) \ dt \ + \alpha \int_0^T \widehat{g}_{ \alpha}^{\prime}(t)\left(g_{\alpha}^{\prime}(t)-\widehat{g}_{ \alpha}^\prime(t)\right) \ dt \geq 0.~	
\end{eqnarray}
Similarly, we replace $ g_{ \alpha}$ with $ \widehat{g}_{ \alpha}(t)$ and $ \bar{g}$ with $ g_{ \alpha}$ in (\ref{546}) we get,
\begin{eqnarray}\label{nb}
	\int_{0}^{T} \left(\widehat{g}_{ \alpha}(t)-g_{ \alpha}(t)\right)\phi(\ell,t;g_{ \alpha}) \ dt \ + \alpha \int_0^T g_{ \alpha}^{\prime}(t)\left(\widehat{g}_{\alpha}^{\prime}(t)-g_{ \alpha}^\prime(t)\right) \ dt \geq 0.~
\end{eqnarray}
We deduce from (\ref{bn}) and (\ref{nb})  that
\begin{eqnarray}\label{mnb1}
	\alpha \int_0^T \left(\widehat{g}_{ \alpha}^{\prime}(t)-g_{ \alpha}^{\prime}(t)\right)^2 dt \leq \int_0^T \left(\widehat{g}_{ \alpha}(t)-g_{ \alpha}(t)\right) \delta \phi(\ell,t)dt,
\end{eqnarray}
where $\delta \phi(\ell,t)=\phi(\ell,t;g_{ \alpha})-\phi(\ell,t;\widehat{g}_{ \alpha})$ is the solution of the adjoint problem (\ref{4}) with data  $ \delta \xi(t)= \delta u(\ell,t)-\delta \nu(t),$ $ \delta u(\ell,t)=u(\ell,t;g_{\alpha})-u(\ell,t;\widehat{g}_{\alpha})$ and $\delta\nu(t)=\nu(t)-\widehat\nu(t).$
Applying H\"older's inequality on the right-hand side of (\ref{mnb1}) and squaring on both sides, we obtain
\begin{eqnarray}\label{549}
	\alpha^2 \Vert \widehat{g}_{ \alpha}^{\prime}-g_{ \alpha}^{\prime}\Vert^4_{L^2(0,T)}\leq \Vert \widehat{g}_{ \alpha}-g_{ \alpha}\Vert^2_{L^2(0,T)}\Vert \delta \phi(\ell,.)\Vert_{L^2(0,T)}^2.
\end{eqnarray}
The inequalities (\ref{1-10}) and (\ref{12}) further lead to the following estimates
\begin{eqnarray}
	\Vert \widehat{g}_{ \alpha}-g_{ \alpha}\Vert^2_{L^2(0,T)}&\leq& T^2 \Vert \widehat{g}_{ \alpha}^{\prime}-g_{ \alpha}^{\prime}\Vert^2_{L^2(0,T)} \label{548}\\
	\Vert \delta \phi(\ell,.)\Vert^2_{L^2(0,T)}&\leq& \frac{T^2\ell^3}{6}\Vert \delta \phi_{xxt}\Vert^2_{L^2(0,T;L^2(0,\ell))} \label{548i}, \nonumber
\end{eqnarray}
so that making use of (\ref{bvc}),  the estimate \eqref{549} becomes as follows
\begin{eqnarray*}\label{540}
		\alpha^2 \Vert \widehat{g}_{ \alpha}^{\prime}-g_{ \alpha}^{\prime}\Vert^2_{L^2(0,T)}&\leq& \frac{T^4\ell^6}{18r_0\kappa_0}\exp(T)(1+T)\Vert \delta \xi^{\prime}\Vert^2_{L^2(0,T)}\\
	&\leq& \frac{T^4\ell^6}{9r_0\kappa_0}\exp(T)(1+T)\Big[\Vert \delta u_t(\ell,.)\Vert^2_{L^2(0,T)} +\Vert \delta \nu^{\prime}\Vert^2_{L^2(0,T)}\Big] .
\end{eqnarray*} 
By using the trace estimate $ \Vert\delta u_{t}(\ell,.)\Vert_{L^2(0,T)}^2\leq \frac{\ell^3}{3}\Vert \delta u_{xxt}\Vert^2_{L^2(0,T;L^2(0,\ell))}$ and appealing to  (\ref{2.98}), which also holds for $ \delta u,$  we get
\begin{eqnarray*}
	\alpha^2 \Vert \widehat{g}_{ \alpha}^{\prime}-g_{ \alpha}^{\prime}\Vert^2_{L^2(0,T)} \leq C_{\alpha}(\kappa_0)\Vert \widehat{g}_{\alpha}^{\prime}-g_{\alpha}^{\prime}\Vert^2_{L^2(0,T)}  + \frac{T^4\ell^6}{9r_0\kappa_0}\exp(T)(1+T)\Vert \delta \nu^{\prime}\Vert^2_{L^2(0,T)},
\end{eqnarray*}
where $ C_{\alpha}(\kappa_0)=\frac{T^4 \ell^{12}}{81\kappa_0^2 r_0^2}\exp(2T)(1+T)^2$ and $\delta g ^{\prime}(t) = g_{\alpha}(t)-\widehat{ g}_{ \alpha}^{\prime}(t).$

Suppose the lower bound of the internal damping coefficient $ \kappa(x)$ be chosen so that $ C_{\alpha}(\kappa_0)=\frac{\alpha^2}{2},$ that is choosing $ \kappa_0$ as in (\ref{643}) and invoking (\ref{548}), we obtain the stability estimate (\ref{646}) with the stability constant $ C_{ST}.$ Hence the proof.
\hfill{$\Box$}

\par We infer from Theorem \ref{mnk} that the stability estimate (\ref{646}) and the lower bound of the internal damping coefficient $ \kappa(x)$ (with units $kgm^3/s$) are valid for all non-negative values of the external damping coefficient $ \mu(x),$ including the critical case,  $ \mu(x)=0 \ kg/ms.$
The following example illustrates the specific values of lower bound for internal damping coefficient corresponding to four different final times $ T$, which are obtained by choosing  $ \ell=.5 \ m,$ and $r_0=1$ (see, \cite{banks1991}).  Using the specific values of $ \kappa_0$, we can analyze the stability constants $ C_{ST}.$
\begin{table}[htb]
	\renewcommand{\arraystretch}{1.85}
	\centering \begin{tabular}{|c|c|c|c|c|c|c|}
		\hline
		$T $&$\alpha$&$\kappa_0=.0025 \left(\frac{T^2 \exp(T)(1+T)}{\alpha}\right)$&$ C_{ST}= 576 \left(\frac{ \kappa_0T^2}{\exp(T)(1+T)}\right)$\tabularnewline
		\hline
		$.1$ & $10^{-3}$ & $.0303$&$.143$\tabularnewline
		\hline
		&$10^{-2}$&$.155$&$9.02$\tabularnewline
		$.5$ &$10^{-3}$ & $1.55$&$90.25$\tabularnewline
		& $ 10^{-4}$&$15.5$&$902.51$ \tabularnewline
		\hline
		$.75$ & $10^{-3}$ & $5.21$&$455.64$\tabularnewline
		\hline
		$1$ & $10^{-3}$ &$ 13.6$&$1440.91$\tabularnewline
		\hline
	\end{tabular}
\end{table}
In a real application, the value of the parameter of regularization ranges between $10^{-2}$ to $10^{-4}.$ The formula (\ref{643}) shows that the lower bound $\kappa_0$ of the internal damping coefficient $\kappa(x)$ is in the order of $\alpha^{-1}$. Hence, a smaller value of $\alpha$ increases the value of the lower bound for the internal damping coefficient, and it drastically increases the stability constant $C_{ST}$. This is evidently clear from the second to fourth rows of the table for the fixed time $ T=.5 \ s.$ Also, one may notice that the increase in final time $ T$ increases the lower bound $\kappa_0$, and the stability constant $C_{ST}$, which indicates that the stability estimates hold only for small intervals of time $T$. The rows corresponding to the fixed $ \alpha =10^{-3}$ show that the  stability constant $C_{ST}$ increases drastically when $ T$ increases. Hence, to get the consistent stability estimate for the IBVP-1, the value of the final time $T>0$ must be small, which is reasonable in terms of applications.

\begin{remark}
The assumption \eqref{643} on the Kelvin-Voigt damping coefficient $\kappa(x)=c_dI(x)$ can be justified by fixing the specific values of various coefficients in \eqref{1-1} and utilizing the estimations of the damping coefficients derived from the dynamic experiments (see, \cite{banks1991})  of real applications.    For a beam of length $l=1,$ moment of inertia $I(x)= 1.64 \times 10^{-9} ,$ mass density $\rho(x)=1.02$ and Young's modulus $E(x)=2.68 \times 10^{10} ,$ the damping coefficients are  estimated  in (\cite{banks1991})  as  $\mu=1.7561 $ and $c_d= 2.05 \times 10^{5} .$  Thus,  for the choice of $\alpha= 10^{-2}$ and $T=0.0133,$ we can see from \eqref{643} that $c_d\geq \kappa_0/I(x)=39597.48.$
\end{remark}	
%\textcolor{blue}{
%\begin{remark}
%The assumption \eqref{643} on the Kelvin-Voigt damping coefficient $\kappa(x)=c_dI(x)$ can be justified by fixing the specific values of various coefficients in \eqref{1-1} and utilizing the estimations of the damping coefficients derived from the dynamic experiments (see, \cite{banks1991})  of real applications.    For a beam of length $l=1\ m,$ moment of inertia $I(x)= 1.64 \times 10^{-9} \ m^4,$ mass density $\rho(x)=1.02\ kg/m$ and modulus $E(x)=2.68 \times 10^{10} \ N/m^2,$ the damping coefficients are  estimated  in (\cite{banks1991})  as  $\mu=1.7561 kg/m s$ and $c_d= 2.05 \times 10^{5} \ kg/m s.$  Thus,  for the choice of $\alpha= 10^{-2}$ and $T=0.0133 s,$ we can see from \eqref{643} that $c_d\geq \kappa_0/I(x)=39597.48.$
%\end{remark}	}
In the case of IBVP-2, we obtain a stability estimate under a sufficient condition on the parameter of regularization $\alpha>0$.
\begin{theorem} \label{st2}
	Assume the conditions given in (\ref{3-5}) and Theorem \ref{t1} hold true. Let $ g_{\alpha}, \ \widehat{g}_{\alpha}\in \mathcal{G}_3$ are unique minimizers of the regularized Tikhonov functional $\mathcal J_{2\alpha}$ defined by (\ref{55}) corresponding to the measured outputs $ \omega, \ \widehat{\omega} \in H^2(0,T),$ respectively. If the parameter of regularization $\alpha$ satisfies the condition
	\begin{eqnarray}\label{345}
		\alpha^2 >  C_{9}^2 C_{10}^2,
	\end{eqnarray}
	where $ C_{9}^2 \ \mbox{and} \ C_{10}^2$ are defined in the proof.
	Then the following stability estimate holds:
	\begin{eqnarray}\label{768}
		\Vert \widehat{g}_{\alpha}-g_{ \alpha}\Vert^2_{L^2(0,T)}&\leq& \widetilde {C}_{ST}\Vert \widehat{\omega}-\omega \Vert^2_{H^2(0,T)},
	\end{eqnarray}
	where $\widetilde{C}_{ST}=\frac{4T^{13}\ell^3 C_6^2 \exp(T/\rho_0)}{ \left(\alpha^2-C_{10}^2 C_{9}^2\right)3r_0}.$
\end{theorem}
{\bf Proof.}
By repeating  the similar steps done in Theorem \ref{mnk} for the variational inequality (\ref{547}), we get
\begin{eqnarray}\label{vi2}
	\alpha\int_0^T \left(\widehat{g}_{\alpha}^{\prime \prime \prime}(t)-g_{\alpha}^{\prime \prime \prime}(t)\right)^2 dt \leq \int_0^T \left(\widehat{g}_{\alpha}(t)-g_{ \alpha}(t)\right) \delta \varphi(\ell,t) dt,
\end{eqnarray}
where $ \delta \varphi(\ell,t)= \varphi(\ell,t;g_{ \alpha})-\varphi(\ell,t;\widehat{g}_{ \alpha})$ is the solution of the adjoint problem (\ref{2}) with $\delta \theta(t)=r(0)\delta u_{xx}(0,t)+\kappa(0)\delta u_{xxt}(0,t)+\delta\omega(t).$ The repeated application of H\"older's inequality gives
\begin{eqnarray} \label{0987}
	\Vert \widehat{g}_{ \alpha}-g_{ \alpha}\Vert^2_{L^2(0,T)}\leq  T^6 \Vert \widehat{g}_{ \alpha}^{\prime \prime\prime}-g_{ \alpha}^{\prime \prime \prime}\Vert^2_{L^2(0,T)}, \  \mbox{for} \  g_{\alpha}, \ \widehat{g}_{\alpha}\in \mathcal{G}_3,
\end{eqnarray}
and using  $ \Vert \delta \varphi(\ell,.)\Vert^2_{L^2(0,T)}\leq \frac{\ell^3}{3}\Vert \delta \varphi_{xx}\Vert^2_{L^2(0,T;L^2(0,\ell))},$ we deduce from \eqref{vi2} that
\begin{eqnarray*}\label{5401}
	\alpha^2 \Vert \widehat{g}_{ \alpha}^{\prime \prime\prime}-g_{ \alpha}^{\prime \prime \prime}\Vert^4_{L^2(0,T)}&\leq& \Vert \widehat{g}_{ \alpha}-g_{ \alpha}\Vert^2_{L^2(0,T)}\Vert \delta \varphi(\ell,.)\Vert_{L^2(0,T)}^2 \nonumber \\  
	&	\leq& \frac{T^6 \ell^3}{3}\Vert \delta \varphi_{xx}\Vert^2_{L^2(0,T;L^2(0,\ell))} \Vert \widehat{g}_{ \alpha}^{\prime \prime\prime}-g_{ \alpha}^{\prime \prime \prime}\Vert^2_{L^2(0,T)} .
\end{eqnarray*}	
In virtue of  the estimates (\ref{097}), (\ref{4.9}) and (\ref{076}), one can get
\begin{eqnarray*}
	\lefteqn{\alpha^2\Vert \widehat{g}_{\alpha}^{\prime \prime \prime}-g_{ \alpha}^{\prime \prime \prime}\Vert^2_{L^2(0,T)} \leq \frac{T^7\ell^3 C_6^2}{3r_0}\exp(T/\rho_0) \left[\Vert \delta \theta ^{\prime}\Vert^2_{L^2(0,T)}+\Vert \delta \theta^{\prime \prime}\Vert^2_{L^2(0,T)}\right]}\nonumber \\
	&\leq& \frac{2T^7\ell^3 C_6^2}{3r_0}\exp(T/\rho_0)\Big[ C_{7}^2 \Big(\Vert \delta u_{tttt}\Vert^2_{L^2(0,T;L^2(0,\ell))}+2\Vert \delta u_{ttt}\Vert^2_{L^2(0,T;L^2(0,\ell))}\nonumber \\ &&+\Vert \delta u_{tt}\Vert^2_{L^2(0,T;L^2(0,\ell))}+\Vert \delta g^{\prime \prime } \Vert^2_{L^2(0,T)}+\Vert \delta g^{\prime}\Vert^2_{L^2(0,T)}\Big) +2\Vert \delta \omega \Vert^2_{H^2(0,T)}\Big],
\end{eqnarray*}
where $\delta g(t) = g_{\alpha}(t)-\widehat{ g}_{ \alpha}(t),$ the constants $ C_6^2= \frac{2\ell^3}{3}\Big( \max (\mu_1^2,\rho_1^2)+\rho_1\max ( 1/T,\,T/3)\Big),$  and $  C_{7}^2=2\ell^2 \max \left( 1,\frac{2\ell}{3}(\rho_1^2+\mu_1^2)\right).$
By employing the regularity estimates (\ref{09}), (\ref{4.7}), (\ref{4.1}) for $ \delta u_{tt},\ \delta u_{ttt}$ and $ \delta u_{tttt}$ respectively,  we deduce that
\begin{eqnarray}\label{569}
	\alpha^2\Vert \widehat{g}_{\alpha}^{\prime \prime \prime}-g_{ \alpha}^{\prime \prime \prime}\Vert^2_{L^2(0,T)}
	&\leq& \frac{2T^7\ell^3 C_6^2}{3r_0}\exp(T/\rho_0) C_{7}^2 \Big[\frac{1}{2\rho_0}\Big(3C_5^2 \exp(C_5^2 T)
	+ C_1^2\ell^3(1+T)\Big)+1\Big]\Vert \delta g \Vert^2_{H^3(0,T)}  \nonumber \\ &&+\frac{4T^7 \ell^3}{3r_0}C_6^2 \exp(T/\rho_0)\Vert \delta \omega \Vert^2_{H^2(0,T)},
\end{eqnarray}
where $ C_5^2$ is the constant defined in Theorem \ref{8}.
For any $ g \in \mathcal{G}_3,$ we obtain that
\begin{eqnarray}\label{567}
	\Vert g \Vert ^2_{H^3(0,T)}\leq C_{9}^2\Vert g^{\prime \prime \prime} \Vert^2_{L^2(0,T)},
\end{eqnarray}
where $ C_9^2=\left( 1+T^2+T^4+T^6\right).$
Substituting (\ref{567}) in (\ref{569}), we arrive at
\begin{eqnarray*}
	\alpha^2\Vert \widehat{g}_{\alpha}^{\prime \prime \prime}-g_{ \alpha}^{\prime \prime \prime}\Vert^2_{L^2(0,T)}\leq C_{10}^2 C_{9}^2\Vert \delta g^{\prime \prime \prime} \Vert^2_{L^2(0,T)} + \frac{4T^7\ell^3}{3r_0}C_6^2 \exp(T/\rho_0)\Vert \delta \omega \Vert^2_{H^2(0,T)},
\end{eqnarray*}
where  $ C_{10}^2 =\frac{2T^7\ell^3 C_6^2}{3r_0}\exp(T/\rho_0) C_{7}^2 \Big[\frac{1}{2\rho_0}\Big(3C_5^2 \exp(C_5^2 T)+ C_1^2 \ell^3(1+T)\Big)+1\Big].$
\\ Choosing $ \alpha^2> C_9^{2}C_{10}^{2},$  %$$ \alpha> \frac{3 T^{11}\ell^3 C_6^2C_{7}^2}{r_0\rho_0^2} \left[C_5^4+C_1^2 (\ell^3+2)\right], $$
we conclude the stability result (\ref{768}) through (\ref{0987}).  \hfill{$\Box$}

The condition (\ref{345}) is simple  to test and implement, and it does not impose a significant constraint. For instance, when $ T=.04 , \ \ell=.4 , \ \rho_1=1 , \ \mu_1=1 , \ r=20 $ and $ \kappa_0=1 ,$ we obtain $ \alpha^2>9.37\times 10^{-9}.$ In  most of the physical experiments, $ \alpha$ varies from $10^{-2} $ to $ 10^{-4},$ and so to validate the condition $\alpha^2>9.37\times 10^{-9},$ one can choose $ \alpha$ as $10^{-2}.$  In this case, we obtain the stability constant $ \widetilde{C}_{ST}=3.27\times 10^{-17},$ which is comparatively small. Therefore, Theorem \ref{st2}   provides a significant stability estimate for the determination of the shear force in terms of a feasible condition on the regularization parameter $\alpha.$ 
\section{Conclusions}
In this work, we studied two inverse problems of identifying unknown transverse shear force in the damped Euler-Bernoulli beam from the boundary data given by measured deflection and bending moment. We analyzed the effect of the Kelvin-Voigt damping  $(\kappa(x)u_{xxt})_{xx}$ and external damping  $\mu(x)u_t$ in the solvability of direct and inverse problems. Though the Kelvin-Voigt damping brings the complicated moment-dependent boundary conditions, the damping effects helped to obtain the enhanced regularity in weak and regular weak solutions of the direct problem with less regular boundary data $g$ in comparison to the undamped Euler-Bernoulli beam.   The same scenario has happened for the inverse problem. We proved the existence of solutions to the inverse problems with less regularity on admissible source inputs $\mathcal G_1$ and  $\mathcal G_2.$ Besides, the Fr\'echet derivative of the Tikhonov functionals is expressed in terms of the adjoint of the direct problems and determined the class of admissible shear forces for which functional's Fr\'echet gradients are Lipschitz continuous. The last result ensures that the sequence $\{\mathcal{J}(g^{(n)})\}$  is monotone,  which assists in devising a gradient-based numerical approach for determining the unknown shear force. A forthcoming paper will carry out the effect of the Kelvin-Voigt damping in the numerical study of reconstructing the shear force in the damped Euler-Bernoulli beam. Finally, we derived local Lipschitz-type stability estimates for the unknown shear force when the Kelvin-Voigt damping coefficient and regularization parameter meet plausible constraints. These findings further show that the stability holds in the absence of external damping effects, and the value of the final time $T>0$ is reasonably small. 

\section{Acknowledgments}
The work of the first author is supported by the National Board for Higher Mathematics, Govt. of India through the research grant No:02011/13/2022/R\&D-II/10206. The research of the second author has been supported by the Scientific and Technological Research Council of Turkey (TUBITAK). The authors would like to thank the reviewer's comments for improving the article further.

\end{document}